\documentclass[letterpaper,11pt]{article}
\usepackage[latin1]{inputenc}
\usepackage{latexsym}
\usepackage[dvips]{graphicx}
\usepackage{amsfonts}
\usepackage{epsfig}
\usepackage{geometry}
\usepackage{amsmath}
\usepackage{amssymb}

\newtheorem{thm}{Theorem}
\newtheorem{prop}[thm]{Proposition}
\newtheorem{lem}[thm]{Lemma}
\newtheorem{cor}[thm]{Corollary}

\newtheorem{ejem}{Example}

\newcommand{\N}{\mathbf{N}}

\newcommand{\R}{\mathbf{R}}

\newcommand{\luno}{\ell^1}
\newcommand{\ltwo}{\ell^2}

\newcommand{\C}{\mathcal{C}}

\renewcommand{\S}{\mathcal{S}}

\newcommand{\f}{\varphi}

\newcommand{\eps}{\varepsilon}

\newcommand{\mles}{{\em M-LES}}
\newcommand{\ces}{{\em CES}}

\newcommand{\limty}[1]{\lim\limits_{#1\to\infty}}
\newcommand{\limn}{\limty{n}}
\def\wto{\rightharpoonup}

\newcommand{\dem}{\bigskip\noindent{\bf Proof. }}
\newcommand{\bx}{\hfill$\blacksquare$\\ \bigskip}
\newcommand{\wbx}{\hfill$\square$\\ \bigskip}
\newcommand{\theend}{\hfill$\square$\\ \bigskip}

\def\interior{\hbox{int}}
\def\argmin{\hbox{Argmin}}

\def\dom{\hbox{dom}}
\newcommand{\sfrac}[2]{\hbox{$\frac{#1}{#2}$}}

\def\Ya{A_{\lambda}}

\def\noi{\noindent}

\begin{document}

\title{Evolution equations for maximal monotone operators: asymptotic analysis in continuous and discrete time}

\author{Juan Peypouquet \\
Departamento de Matem\'atica, Universidad T\'ecnica Federico Santa Mar\'\i a\\
Av. Espa\~na 1680, Valpara\'\i so, Chile\\
{\small juan.peypouquet@usm.cl}\\
Sylvain Sorin \\
Equipe Combinatoire et Optimisation, CNRS FRE 3232,  Facult\'e de Math\'ematiques,\\
Universit\'e P. et M. Curie - Paris 6,
  175 Rue du Chevaleret,
75013 Paris\\  and Laboratoire d'Econom{\'e}trie, Ecole Polytechnique, France \\
{\small  sorin@poly.polytechnique.fr}}

\maketitle

\small

\begin{abstract}
This survey is devoted to the asymptotic behavior of solutions of evolution equations generated by maximal monotone operators in Hilbert spaces. The emphasis is in the comparison of  the continuous time trajectories to sequences generated by implicit or explicit discrete time schemes.
The analysis covers weak convergence for the average process, for the process itself  and strong convergence and aims at highlighting the main ideas and unifying the proofs. We further make the connection with the analysis in terms of almost orbits that allows for a broader scope.
\end{abstract}

\tableofcontents

\section*{Introduction}

Discrete and continuous dynamical systems governed by maximal monotone operators have a great number of applications in
optimization, equilibrium, fixed-point theory, partial differential equations, among others.\\

We are specially concerned about the connection between  continuous and  discrete models. This connection occurs at two
levels:
\begin{enumerate}
    \item On  a compact interval, one  approximates a  continuous-time trajectories by interpolation of some sequences computed via discretization.   By  considering  vanishing step size   this construction is used to prove existence results and to approximate the trajectories numerically.
    \item  Another approximation is in the long term,  were we compare  asymptotic properties of a continuous trajectory to similar    asymptotic properties of a given   path defined  inductively trough a sequence of values and step sizes.
\end{enumerate}
It is important to mention that  some  estimations (eg. Kobayashi) can be useful for both purposes.\\

The literature on this subject is huge   but  lot of the arguments turn out
to be pretty much the same. Therefore, we intend to give a concise yet complete compendium of the results available, with an emphasis on the techniques and the way the enter in the proofs.\\
Most of the properties will be established  in the framework of Hilbert spaces since  our aim is to emphasize  unity in terms of tools and approach. A lot of results  can be extended  but in most of  the case under specific assumptions.  With no aim for completeness, we have included several
references to the corresponding results in Banach spaces that we think might be useful.\\

The paper is organized as follows: In section 1 we recall the basic properties of maximal monotone operators along with
some examples. Section 2 deals with the associated dynamic approach. We present the existence results for the
differential inclusion $\dot u\in -Au$ and global properties of implicit and explicit discretizations. Section 3
establishes the convergence of the value $f(u)$ in the case of an operator of the form $A=\partial f$. In section 4 we
describe  general results on weak convergence: tools, arguments, characterization of the weak limits. Section 5 is
devoted to weak convergence in average and Section 6 is concerned with weak convergence, especially for demipositive
operators. In section 7 we present the, mostly geometric, conditions ensuring that the convergence is strong. Section 8 deals with
asymptotic equivalence  and explains some apparently hidden relationships between certain continuous- and discrete-time
dynamical systems. Finally, section 9 contains some concluding remarks.\\

\section{Preliminaries}

The purpose of this section is to introduce notations and to recall basic results.

\subsection{Maximal monotone operators}

Let $H$ be a real Hilbert space with inner product $\langle\cdot,\cdot\rangle$ and norm $\|\cdot\|$. An {\em operator}
is a set-valued mapping $A:H\rightrightarrows H$ whose domain $$D(A)=\{u\in H: Au\neq \emptyset\}$$ is nonempty. For
convenience of notation, sometimes we will identify $A$ with its graph by writing $[u,u^*]\in A$ for $u^*\in Au$. The
operator $A^{-1}$ is defined by its graph: $[u,u^*]\in A^{-1}$ if, and only if, $[u^*,u]\in A$. \\
An operator $A:H
\rightrightarrows H$ is {\em monotone} if one has
\begin{equation}\label{Inec:monotone}
\langle x^*-y^*,x-y\rangle\ge 0
\end{equation}
for all $[x,x^*], [y,y^*]\in A$.\\
 A monotone operator is {\em maximal} if its graph is not properly contained in the
graph of any other monotone operator.  Observe that if $A$ is monotone (resp. maximal monotone)  then so are $A^{-1}$
and $\lambda A$ if $\lambda>0$.

\begin{lem}\label{lem:mm_sufficient}
Let $A$ be a maximal monotone operator. A point $[x,x^*]\in H\times H$ belongs to the graph of $A$ if, and only if,
$$\langle x^*-u^*,x-u\rangle\ge 0\qquad\hbox{for all}\quad [u,u^*]\in A.$$
\end{lem}

\dem If $[x,x^*]\in A$ the inequality holds by monotonicity. Conversely, if $[x,x^*]\notin A$, then the set
$A\cup\{[x,x^*]\}$ is the graph of a monotone operator that extends $A$, which contradicts maximality. \bx

An operator $A:H
\rightrightarrows H$ is {\em nonexpansive} if one has
\begin{equation}\label{nonex}
\| x^*-y^*\|\leq \| x-y\|
\end{equation}
for all $[x,x^*], [y,y^*]\in A$. Observe that a nonexpansive operator is single-valued on its domain. \\

Let $I$ be the identity mapping on $H$. For $\lambda>0$, the {\em resolvent} of $A$ is the operator
$$J^A_\lambda=(I+\lambda A)^{-1}.$$

\begin{thm}\label{Thm:Minty}
Let $A:X\rightrightarrows X$. Then
\begin{itemize}
    \item [i)] $A$ is monotone if, and only if, $J^A_\lambda$ is  nonexpansive for each $\lambda>0$.
    \item [ii)] A monotone operator $A$ is maximal if, and only if, $I+\lambda A$ is surjective for each $\lambda>0$.
\end{itemize}
\end{thm}

\dem \ \\
$i)$ Let $A$ be monotone,  $[x,x^*], [y,y^*]\in A$ and $\lambda>0$. \\
Inequality  (\ref{Inec:monotone}) implies
\begin{equation}\label{incnorm}
\| x- y \| \le \| x-y + \lambda ( x^*- y ^*) \|,  \quad  \forall \lambda \geq 0
\end{equation}
which is the non expansiveness of $J^A_\lambda$.\\
Conversely, (\ref{incnorm}) leads to
$$
\lambda \langle x^*-y^*,x-y\rangle + \lambda^2 \|x^*- y ^*\|^2 \geq 0
$$
hence implies  (\ref{Inec:monotone}) by dividing by $\lambda$ and  letting $\lambda \to 0$.

\noi  $ii)$ It is enough    to prove the result for $\lambda=1$. Given $z_0\in H$, we will
find $x_0\in H$ such that $\langle
    y-(z_0-x_0),x-x_0\rangle\geq 0$ for all $[x,y]\in A$ so that maximality of $A$ implies
    $z_0-x_0\in Ax_0$. For ${[x,y]\in A}$, define the  weakly compact
set $C_{x,y}$  by
$$
C_{x,y}=\{x_0\in H:\langle y+x_0-z_0,x-x_0\rangle\geq 0\}.
$$
It suffices to show that the family $\{C_{x,y}\}_{[x,y]\in A}$ has the finite intersection property. To this end take $[x_i,y_i]\in A$ for
$i=1,\ldots,n$. Let $\Delta=\{(\lambda_1,\ldots,\lambda_n):\lambda_i\geq 0; \sum_{i=1}^n\lambda_i=1\}$ denote the
$n$-dimensional simplex and consider the function $f:\Delta\times\Delta\to\R$ given by
$$f(\lambda,\mu)=\mbox{$\sum_{i=1}^n$}\,\mu_i\langle y_i+x(\lambda)-z_0,x(\lambda)-x_i\rangle$$
with $x(\lambda)\!=\!\sum_{i=1}^n\lambda_i x_i$. Clearly $f(\cdot,\mu)$ is convex and continuous while
$f(\lambda,\cdot)$ is linear. The Min-max Theorem (see, for instance, Theorem 1.1 in \cite[Br\'ezis]{Bre}) implies the
existence of $\lambda_0\!\in\!\Delta$ such that
$$\max_{\mu\in\Delta}f(\lambda_0,\mu)=
\max_{\mu\in\Delta}\min_{\lambda\in\Delta}f(\lambda,\mu) \leq\max_{\mu\in\Delta}f(\mu,\mu).$$ \noi Now
monotonicity of $A$ implies
\begin{eqnarray*}
f(\mu,\mu)&=&\hbox{$\sum_{i=1}^n$}\,\mu_i\langle y_i,x(\mu)-x_i\rangle+
\langle x(\mu)-z_0,x(\mu)-x(\mu)\rangle\\
&=&\hbox{$\sum_{i,j=1}^n$}\,\mu_i\mu_j\langle y_i,x_j-x_i\rangle\\
&=&\hbox{$\frac{1}{ 2}\sum_{i,j=1}^n$}\,\mu_i\mu_j\langle y_i-y_j,x_j-x_i\rangle\le  0
\end{eqnarray*}
so that $f(\lambda_0,\mu)\!\leq\! 0$ for all $\mu\!\in\!\Delta$, and taking for $\mu$ the extreme points we get
$\langle y_i\!+\!x(\lambda_0)\!-\!z_0,x(\lambda_0)\!-\!x_i\rangle\!\leq\! 0$ for all $i$, which  is
$x(\lambda_0)\!\in\! \bigcap_{i=1}^nC_{x_i,y_i}$.

Conversely, take $[u,u^*]\in H\times H$ such that $\langle u^*-v^*,u-v\rangle\ge 0$ for all $[v,v^*]\in A$. We shall
prove that $[u,u^*]\in A$. Since $I+A$ is surjective, there is $[\overline{v},\overline{v}^*]\in A$ such that
$\overline{v}+\overline{v}^*=u+u^*$. Then $ \langle u^* - \overline{v}^*, u - \overline{v} \rangle  = -  \|u-\overline{v}\|^2\geq 0$ which implies  $u=\overline{v}$,
$u^*=\overline{v}^*$ and $[u,u^*]\in A$. \bx

{\it Comments }\\
The study of monotone operators started in \cite[Minty]{Min}. See also \cite[Kato]{Kat} for part $i)$ in Banach spaces.
The {\em if} part in $ii)$ holds in Banach spaces, but not the {\em only if} part (see \cite[Hirsh]{Hir}). \theend

\subsection{Examples and properties}

\begin{ejem}{\em Let $C\subset H$ and let $T:C\to H$ be nonexpansive. The operator $A=I-T$ is monotone because
\begin{eqnarray*}
\langle Ax-Ay,x-y\rangle & = & \|x-y\|^2-\langle Tx-Ty,x-y\rangle\\
& \ge & \|x-y\|\left[\|x-y\|\frac{}{}\!-\|Tx-Ty\|\right]\\
& \ge & 0.
\end{eqnarray*}
Maximality depends on whether $T$ can be extended to a nonexpansive function on a set that contains $C$ properly (for
example if $C$ is closed and convex).}\theend
\end{ejem}

\begin{ejem}{\em Let $\Gamma_0(H)$ denote the set of all proper, lower-semicontinuous convex functions $f:H\to\R\cup\{+ \infty\}$.
For $f\in\Gamma_0(H)$, the {\em subdifferential of} $f$ is the operator $\partial f :H\rightrightarrows H$ defined by
$$\partial f(x)=\{x^*\in H:f(z)\ge f(x)+\langle x^*,z-x\rangle\hbox{ for all }z\in H\}.$$
To see that it is monotone, take $x^*\in\partial f(x)$ and $y^*\in\partial f(y)$. Thus
\begin{eqnarray*}
f(y) & \ge & f(x)+\langle x^*,y-x\rangle\\
f(x) & \ge & f(y)+\langle y^*,x-y\rangle.
\end{eqnarray*}
and adding these two inequalities we obtain  $\langle x^*-y^*,x-y\rangle\ge 0$. For maximality, according to Theorem \ref{Thm:Minty}
it suffices to prove that for each $y\in H$ and each $\lambda>0$
there is $x_\lambda\in D(\partial f)$ such that $y\in x_\lambda+\lambda\partial f(x_\lambda)$. Indeed, consider the
{\em Moreau-Yosida approximation} of $f$ at $y$, which is the function $f_\lambda$ defined by
\begin{equation}\label{yos}
f_\lambda(x)=f(x)+\frac{1}{2\lambda}\|x-y\|^2.
\end{equation}
It is proper, lower-semicontinuous, strongly convex and coercive (due to the quadratic term and the fact that $f$ has a
affine minorant). Its unique minimizer $x_\lambda$ satisfies
$$0\in\partial f_\lambda(x_\lambda)=\partial f(x_\lambda)+\frac{1}{\lambda}(x_\lambda-y).$$
That is, $y\in x_\lambda+\lambda\partial f(x_\lambda)$.}\theend
\end{ejem}

The {\em solution set} of $A$ is $\S=A^{-1}0 = \{x \in H;   0  \in Ax\}$. This set is relevant in optimization and
fixed-point theory:
\begin{itemize}
    \item If $A=I-T$, where $T$ is a nonexpansive mapping, then $\S$ is the set of fixed points of $T$.
    \item If $A=\partial f$, where $f$ is a proper lower-semicontinuous convex function then $\S$ is the set of minimizers
        of $f$.
\end{itemize}

Let us describe some topological consequences of maximal monotonicity.

\begin{prop}\label{prop:ws_sw}
Let $A$ be a maximal monotone operator. Then $A$ is sequentially weak-strong and strong-weak closed.
\end{prop}

\dem Take sequences $\{x_n\}$ and $\{x^*_n\}$ in $H$ such that $[x_n,x^*_n]\in A$ for each $n\in\N$ and suppose  that
$x_n\to x$ and $x^*_n\wto x^*$, as $n\to\infty$ (consider $A^{-1}$ for the other case). To prove  that $[x,x^*]\in A$,
recall that by monotonicity,   for all $[u,u^*]\in A$ and all $n\in\N$,  $\langle x_n^*-u^*,x_n-u\rangle\ge 0$. Letting
$n\to\infty$ the convergence assumptions imply that $\langle x^*-u^*,x-u\rangle\ge 0$ for all $[u,u^*]\in A$. Hence
$[x,x^*]\in A$ by Lemma \ref{lem:mm_sufficient}. \bx

\begin{cor}
Let $A$ be maximal monotone. For each $x\in D(A)$ the set $Ax$ is closed and convex. In particular, $\S$ is closed and
convex.
\end{cor}

\dem Proposition \ref{prop:ws_sw} implies $Ax$ is closed for each $x\in D(A)$. To see that $Ax$ is convex, take
$x^*,y^*\in Ax$, $[u,u^*]\in A$ and $\lambda\in(0,1)$. Then $\langle \lambda x^*+(1-\lambda)y^*-u^*,x-u\rangle =
\lambda \langle x^*-u^*,x-u\rangle+(1-\lambda)\langle y^*-u^*,x-u\rangle\ge 0.$ As before, we conclude that $\lambda
x^*+(1-\lambda)y^*\in Ax$ by Lemma \ref{lem:mm_sufficient}. Finally, since $A^{-1}$ is maximal monotone and
$\S=A^{-1}0$, the set $\S$ is closed and convex. \bx


\section{Dynamic approach}\label{section:dynamic-approach}

The following sections address, among others,  the issue of finding zeroes of  a maximal monotone operator $A$. The
strategy is the following: we shall consider some continuous and discrete dynamical systems whose trajectories  may
converge, in some sense and under some conditions,  to points in $\S = A^{-1} 0$. In this section we present these
systems along with some relevant properties.\\

{\bf From now on we assume that $A$ is a maximal monotone operator. }

\subsection{Differential inclusion}

In this section we consider the following differential inclusion:
\begin{equation}\label{continuous}
\left\{\begin{array}{rcccl}
\dot u(t)  & \in  & - Au(t) & \hbox{a.e. on $(0,\infty)$}\\
u(0)=x & \in& D(A).
\end{array}\right.
\end{equation}

A {\em solution} of \eqref{continuous} is an  absolutely continuous function  $u$ from $\R^+$ to $H$ satisfying these
two conditions.\\
Monotonicity implies the following dissipative property:

\begin{lem}\label{lem:basic_continuous_1}
Let $u_1$ and $u_2$ be absolutely continuous functions satisfying $\dot u_i(t)\in -Au_i(t)$ almost everywhere on
$(0,T)$. Then the function $t\mapsto\|u_1(t)-u_2(t)\|$ is decreasing on $(0,T)$.
\end{lem}

\dem For $t\in (0,T)$ define $\theta(t)=\frac{1}{2}\|u_1(t)-u_2(t)\|^2$. The hypotheses give $\dot\theta(t)=\langle\dot
u_1(t)-\dot u_2(t),u_1(t)-u_2(t)\rangle\le 0$ for almost every $t$. \bx

Immediate consequences are the following:

\begin{cor}\label{cor:basic_continuous_2}
Let $y\in\S$ and $u$ be a solution of \eqref{continuous}.
Then  $\limty{t}\|u(t)-y\|$ exists.
\end{cor}

\begin{cor}
There is at most one solution of \eqref{continuous}.
\end{cor}

Another aspect of dissipativity is the following:

\begin{prop}\label{u'dec}
$ \| \dot  u (t) \|$  is decreasing.
\end{prop}

\dem Lemma \ref{lem:basic_continuous_1} implies  that for any $h>0$ and  $s<t$
$$
\| u(t+h ) - u(t) \| \leq  \| u(s+h ) - u(s) \|
$$
hence the result by dividing by $h$ and taking the limit as $h \to 0$.\bx

We shall present two approaches for  the existence of a solution  of \eqref{continuous}. The first one uses the {\it Yosida approximation} and is the best-known in the theory of
optimization in Hilbert spaces. The second one uses {\it proximal sequences} to approximate the function $u$. It is popular
in the field of partial differential equations since it works naturally in arbitrary Banach spaces.\\

But before doing so, and assuming for a moment that the differential inclusion \eqref{continuous} {\em does} have a
solution, observe that by Lemma \ref{lem:basic_continuous_1}, for each $t\geq0$ the mapping $x \mapsto u(t)$ defines a
non expansive function from $D(A)$ to itself that can be continuously extended to a map $S_t$ from $\overline{D(A)}$ to
itself. The family $\{S_t\}_{t\geq 0}$ is the  {\it semi-group}  generated by $A$ and satisfies:
\begin{itemize}
    \item [i)] $S_0 = I$ and $S_t\circ S_r = S_{t+r}$;
    \item [ii)] $\| S_tx- S_ty \| \leq \| x-y\|$;
    \item [iii)] $\lim\limits_{t\to 0} \| x - S_t x \| =  0$.
\end{itemize}
Reciprocally, given a {\it continuous semi-group of contractions}  i.e. satisfying i), ii) and  iii),  from a closed convex subset $C$ to itself,  there exists a  {\it generator},
namely a maximal operator $A$ with $C = \overline {D(A)}$ such that $S_tx$ coincides with $u(t)$  for $x\in D(A)$, see
\cite[Br\'ezis]{Bre}.\\

{\bf We will use hereafter both notations $u(t)$ and $S_tx$.}

\subsection{Approach through the Yosida approximation.}

\subsubsection{The Yosida Approximation}

\noi  Recall that the resolvent is $J^A_{\lambda}$. The {\em Yosida approximation} of
 $A$ is the single-valued maximal monotone operator
$A_{\lambda}$, $\lambda >0$, defined by
$$A_{\lambda} = \frac{1}{\lambda}(I- J^A_{\lambda}).$$
Since $J^A_{\lambda}$ is nonexpansive and everywhere defined, $A_\lambda$ is monotone (see example 1 above) and maximal
(using Lemma \ref{lem:mm_sufficient}). It is also clear that $A_\lambda$ is Lipschitz-continuous with
constant $2/\lambda$. Observe that $\S=A^{-1}0 = A_{\lambda}^{-1}0 $ for all $\lambda >0$.\\

\noi For a closed convex set $C\subset H$ and a point $x\in H$ we denote by $P_Cx$ the orthogonal projection of $x$
onto $C$. The {\em minimal section} of $A$ is the operator $A^0$ defined by $A^0x=P_{Ax}0$, which is clearly monotone but
not necessarily maximal.\\

\noi The following results summarize the main properties of the resolvent and the Yosida approximation. They can be
found in \cite[Br\'ezis]{Bre} (see also \cite[Barbu]{Bar} for Banach spaces).

\begin{prop} With the notation introduced above we have the following:
\begin{enumerate}
    \item $A_\lambda x\in AJ^A_\lambda x$
    \item $\|A_\lambda x\|\le \|A^0x\|$, $\|A_\lambda x\|$ is nonincreasing in $\lambda$ and $\lim\limits_{\lambda\to 0}\|A_\lambda x\|\to\|A^0x\|$.
    \item $\lim\limits_{\lambda\to 0}J^A_\lambda x=x$.
    \item If $x_\lambda\to x$ and $A_\lambda x_\lambda$ remains bounded as $\lambda\to 0$, then $x\in D(A)$.
Moreover, if $y$ is a cluster point of $A_\lambda x_\lambda$  as $\lambda\to 0$, then $y\in Ax$.
    \item $A^0$ characterizes $A$ in the following sense: If $A$ and $B$ are maximal monotone with common
domain and $A^0=B^0$, then $A=B$.
    \item $\lim\limits_{\lambda\to 0}A_\lambda x= A^0x$ and $\overline{D(A)}$, the (strong) closure of $D(A)$, is convex.
\end{enumerate}
\end{prop}

\subsubsection{The existence result}

The main result is the following:

\begin{thm}\label{exuni}
There exists a unique absolutely continuous function $u:[0, +\infty)\to H$ satisfying \eqref{continuous}. Moreover,
\begin{enumerate}
    \item $\dot u\in L^\infty(0,\infty;H)$ with $\|\dot{u}(t)\|\le\|A^0x\|$ almost everywhere.
    \item $u(t)\in D(A)$ for all $t\ge 0$ and $\|A^0u(t)\|$ decreases.
    \item $A^0u(t)$ is continuous from the right and $u(t)$ admits a right-hand derivative for all $t\ge 0$; namely
        $\dot{u}(t^+)=-A^0u(t)$ (lazy behavior)
\end{enumerate}
\end{thm}

\noi The problem of finding a trajectory satisfying \eqref{continuous} was first posed and studied in
\cite[Komura]{Kom} and \cite[Crandall and Pazy]{CrP2}. The classical proof can be found in \cite[Br\'ezis]{Bre}. The
idea is to consider the differential inclusion \eqref{continuous} with $A=\Ya$, which has a solution $u_{\lambda}$ by
virtue of the Cauchy-Lipschitz-Picard Theorem. Then one proves  first that, as $\lambda\to 0$, $u_{\lambda}$ converges
uniformly on compact intervals to some $u$, then that $u$  satisfies \eqref{continuous} for the original $A$. The
following estimation plays a crucial role in the proof and is interesting on its own:
\begin{equation}\label{eqBre}
\| u_{\lambda} (t) -  u(t)\| \leq  2\| A^0 (u_0)\| \sqrt { \lambda t}.
\end{equation}
Finally $u$ is proved to have the properties enumerated in Theorem \ref{exuni}.\\

{\it Comments}\\
The same method can be extended to Banach spaces $X$ such that $X$ and $X^*$ are uniformly convex (see
\cite[Kato]{Kat}).\theend

\subsection{Approach through proximal sequences.}

\subsubsection{Proximal sequences}

Let $\{\lambda_n\}$ be a sequence of positive numbers or {\em stepsizes}. $\{x_n\}$ is a {\em proximal sequence} if it
satisfies
\begin{equation}\label{prox}
\left\{\begin{array}{rcl}
\displaystyle \frac{x_n-x_{n-1}}{\lambda_n}  & \in  & - Ax_n\qquad\hbox{for all $n\ge 1$}\\
x_0 & \in & D(A).
\end{array}\right.
\end{equation}
In other words,
\begin{equation}\label{proxex}
x_n=(I+\lambda_n A)^{-1}x_{n-1}=J_{\lambda_n}^Ax_{n-1}.
\end{equation}
The existence of such a sequence follows from Theorem \ref{Thm:Minty}. Observe that the first inclusion in \eqref{prox}
can be seen as an implicit discretization of the differential inclusion \eqref{continuous}, called also a backward
scheme. The {\em velocity} at stage $n$ is $$y_n=\displaystyle \frac{x_n-x_{n-1}}{\lambda_n}.$$

{\it Comments}\\
\noi The notion of proximal sequences and the term {\it proximal} were introduced in \cite[Moreau]{Mor} for $A=\partial
f$. In that case, finding $x_n$ corresponds to minimizing the Moreau-Yosida approximation of $f$ at $x_{n-1}$ (see
\eqref{yos}), namely
$$f_{\lambda_n}(x)=f(x) + \frac{1}{2\lambda _n}\|x- x_{n-1}\|^2.$$\theend

Monotonicity implies the following properties:
\begin{lem}\label{lem:basic_prox_1} The sequence $\|y_n\| $ is decreasing.
\end{lem}

\dem The inequality $\langle y_n - y_{n-1}, x_n - x_{n-1} \rangle \le 0$ implies $ \langle y_n - y_{n-1}, y_n\rangle
\le 0$ and therefore $\|y_n\|\le\|y_{n-1}\|$. \bx

This is the counterpart of $\|\dot u (t) \|$ decreasing, Proposition \ref{u'dec}.

\begin{lem}\label{lem:basic_prox_2} Let  $x\in \S$. Then
$\|x_n-x\|^2 + \lambda_n^2 \| y_n\|^ 2 \leq \| x_{n-1} - x \|^2$.
\end{lem}

\dem Simply observe that
\begin{eqnarray*}
\|x_{n-1}-x\|^2 & = & \|x_{n-1}-x_n\|^2+\|x_n-x\|^2+2\langle x_{n-1}-x_n,x_n-x\rangle\\
& \ge & \lambda_n^2\|y_n\|^2+\|x_n-x\|^2
\end{eqnarray*}
since $\langle x_{n-1}-x_n,x_n-x\rangle \ge 0 $ by monotonicity when $x$ is in $\S$.\bx

An immediate consequence is the following:

\begin{cor}\label{lem:basic_prox_3}
Let $x \in \S$. The sequence $\|x_n  - x\|^2$ is decreasing, thus convergent.
\end{cor}

Notice the similarity with Corollary \ref{cor:basic_continuous_2}.\\

\subsubsection{Kobayashi inequality} The following inequality, due to  Kobayashi \cite{Kob}, provides an
estimation for the distance between two proximal sequences
$\{x_k\}$ and $\{\widehat x_l\}$,  with
stepsizes $\{\lambda_k\}$ and $\{\widehat\lambda_l\}$, respectively. \\
{\bf We use the following notation throughout the
paper}:
$$\sigma_k=\sum_{i=1}^k\lambda_i\qquad \mbox{and} \qquad  \tau_k=\sum_{i=1}^k\lambda_i^2$$ (similarily for $\widehat\sigma_l$ and
$\widehat\tau_l$).

\begin{prop}[Kobayashi inequality]\label{Kobayashi}
Let
 $\{x_k\}$ and $\{\widehat x_l\}$ be two  proximal sequences.  If $u\in D(A)$, then
\begin{equation}\label{intro2_eveq}
\|x_k-\widehat x_l\|\le
    \|x_0-u\|+\|\widehat x_0-u\|+
    \|A^0u\|\sqrt{(\sigma_k-\widehat\sigma_l)^2+\tau_k+\widehat\tau_l}.
\end{equation}
\end{prop}

\noi We first prove the following auxiliary result:

\begin{lem}\label{kobalemma}
Let  $[u_1,\ v_1],\ [u_2,\ v_2]\in A$ and $\lambda,\ \mu>0$, then
$$(\lambda+\mu)\|u_1-u_2\|\le \lambda\|u_2+\mu v_2-u_1\| + \mu\|u_1+\lambda v_1-u_2\|.$$
\end{lem}

\dem Write $u=u_1-u_2$. Then
\begin{eqnarray*}
(\lambda+\mu)\|u_1-u_2\|^2 & = & \lambda\langle u_2-u_1, -u\rangle + \mu\langle u_1-u_2,u\rangle\\
& = & \lambda\langle u_2+\mu v_2-u_1,-u\rangle + \mu\langle u_1+\lambda v_1-u_2, u\rangle\\
& & \hskip20pt     + \lambda\mu\langle v_2-v_1, u_1-u_2\rangle\\
& \le & \left[\lambda\|u_2+\mu v_2-u_1\| + \mu\|u_1+\lambda v_1-u_2\|\right]\|u_1-u_2\|
\end{eqnarray*}
by monotonicity.\bx

\noindent{\bf Proof of Proposition \ref{Kobayashi}:} To simplify notation set
$$c_{k,l}=\sqrt{(\sigma_k-\widehat\sigma_l)^2+\tau_k+\widehat\tau_l}.$$
The proof will use induction on the pair $(k,l)$.\\
First, let us establish inequality \eqref{intro2_eveq} for the pair $(k,0)$ with $k\ge 0$. Monotonicity implies, using (\ref{incnorm})
 that
$$
\| x_1- u \|  \le \| x_0 - u - \lambda_1  A ^{0} u \|
$$
and
$$
\| x_1- u \|  \le \| x_0 - u\|  + \lambda_1 \| A ^{0} u \|.
$$
Inductively we obtain
$$
\| x_k- u \|  \le \| x_0 - u\|  + \sigma_k \| A ^{0} u \|.
$$
thus
\begin{eqnarray*}
\|x_k-\widehat x_0\| & \le & \|x_k-u\|+\|u-\widehat x_0\|\\
& \le & \|x_0-u\|+\sigma_k\|A^0u\|+\|\widehat x_0-u\|\\
& \le & \|x_0-u\|+\|\widehat x_0-u\|+c_{k,0}\|A^0u\|
\end{eqnarray*}
because $\sigma_k\le c_{k,0}$. In a similar fashion we prove the inequality for $(0,l)$ with $l\ge 0$.\\

Now suppose \eqref{intro2_eveq} holds for $(k-1,l)$ and $(k,l-1)$. According to Lemma \ref{kobalemma},
$$(\lambda_k+\widehat\lambda_l)\|x_k-\widehat x_l\| \le \lambda_k\|\widehat
x_l+\widehat\lambda_l\widehat y_l-x_k\|+\widehat\lambda_l\|x_k+\lambda_ky_k-\widehat x_l\|.$$ Setting
$\alpha_{k,l}=\displaystyle \frac{\widehat\lambda_l}{\lambda_k+\widehat\lambda_l}$ and $\beta_{k,l}= 1-\alpha_{k,l}=
\displaystyle \frac{\lambda_k}{\lambda_k+\widehat\lambda_l}$ we have
\begin{eqnarray}\label{kobaux1}
\|x_k-\widehat x_l\| & \le & \alpha_{k,l}\|x_{k-1}-\widehat x_l\|+\beta_{k,l}\|\widehat x_{l-1}-x_k\|\nonumber\\
& \le & \alpha_{k,l}\left[\|x_0-u\|+\|\widehat x_0-u\|+c_{k-1,l}\|A^0u\|\right]\nonumber\\
& {} & \hskip10pt + \beta_{k,l}\left[\|x_0-u\|+\|\widehat x_0-u\|+c_{k,l-1}\|A^0u\|\right]\nonumber\\
& = & \|x_0-u\|+\|\widehat x_0-u\|+\left[\alpha_{k,l}c_{k-1,l}+\beta_{k,l}c_{k,l-1}\right]\|A^0u\|.
\end{eqnarray}
It only remains to verify that
\begin{equation}\label{kobaux2}
\alpha_{k,l}c_{k-1,l}+\beta_{k,l}c_{k,l-1}\le c_{k,l}.
\end{equation}

\noi Cauchy-Schwartz Inequality implies
\begin{eqnarray*}
\alpha_{k,l}c_{k-1,l}+\beta_{k,l}c_{k,l-1} & = &
\alpha^{1/2}_{k,l}(\alpha^{1/2}_{k,l}c_{k-1,l})+\beta^{1/2}_{k,l}(\beta^{1/2}_{k,l}c_{k,l-1})\\
& \le &
(\alpha_{k,l}+\beta_{k,l})^{1/2}(\alpha_{k,l}c^2_{k-1,l}+\beta_{k,l}c^2_{k,l-1})^{1/2}\\
& = & (\alpha_{k,l}c^2_{k-1,l}+\beta_{k,l}c^2_{k,l-1})^{1/2}.
\end{eqnarray*}
On the other hand, notice that $c^2_{k-1,l}=c^2_{k,l}-2\lambda_k(\sigma_k-\widehat\sigma_l)$, while
$c^2_{k,l-1}=c^2_{k,l}+2\widehat\lambda_l(\sigma_k-\widehat\sigma_l)$. Hence,
\begin{eqnarray*}
(\alpha_{k,l}c_{k-1,l}+\beta_{k,l}c_{k,l-1})^2 & \le &
\alpha_{k,l}c^2_{k-1,l}+\beta_{k,l}c^2_{k,l-1}\\
& = & \alpha_{k,l}c^2_{k,l}+\beta_{k,l}c^2_{k,l} -
2(\alpha_{k,l}\lambda_k-\beta_{k,l}\widehat\lambda_l)(\sigma_k-\widehat\sigma_l)\\
& = & c^2_{k,l}.
\end{eqnarray*}
Inequalities \eqref{kobaux1} and \eqref{kobaux2} give \eqref{intro2_eveq}.\bx

{\it Comments}\\
\noi Kobayashi's original inequality also accounts for possible errors in the determination of the proximal sequence,
see \cite{Kob}. Nonautonomous version of the inequality can be found in \cite[Kobayasi, Kobayashi and Oharu]{KKO}
\cite[Alvarez and Peypouquet]{AlP1}.\theend

\subsubsection{The existence result} In general Banach spaces, existence and uniqueness can also be derived by the method in \cite[Crandall and
Liggett]{CrL}, based on the resolvent, which we now present:\\

Set $t\in [0,T]$, $m\in\N$ and consider a proximal sequence with constant stepsizes $\lambda_k\equiv t/m$. The $m$-th
iteration defines a function
$$
u_m(t)=\left(I+\frac{t}{m}A\right)^{-m}x.
$$
 Repeat the procedure for each $m$
to obtain a sequence $\{u_m(t)\}$ of functions from $[0,T]$ to $H$. The following result was proved  in \cite[Crandall
and Liggett]{CrL}:

\begin{thm}\label{thm:proxdiff}
The sequence $\{u_m(t)\}$ defined above converges to some $u(t)$ uniformly on every compact interval $[0,T]$. Moreover,
the function $t\mapsto u(t)$ satisfies \eqref{continuous}.
\end{thm}
\dem Instead of the original proof we present an easier one using Kobayashi's inequality
\eqref{intro2_eveq}\footnote{In fact, Kobayashi's proof is based on a simplification of Crandall and Liggett's
method.}. Fix $N,M\in\N$ and $t,s\in[0,T]$ with $T>0$. Consider two proximal sequences with $\lambda_k= t/N$ and
$\widehat\lambda_l= s/M$ for all $k, l$. Initialize $x_k$ and $\widehat x_l$ both at $x$. Note that  $x_N = u_N(t)$ and
$\widehat x_M = u_M(s)$ hence
$$
\|u_N(t)-u_M(s)\|\le \|A^0x\|\sqrt{(t-s)^2+\sfrac{T^2}{N}+\sfrac{T^2}{M}}.
$$
Thus the sequence $\{u_n\}$ converges uniformly on $[0,T]$ to a function $u$, which is uniformly
Lipschitz-continuous with constant $\|A^0x\|$.

\noi In order to prove that the function $u$ satisfies \eqref{continuous} it suffices to verify that it is an {\em
integral solution} in the sense of B\'enilan \cite{Ben}, which means that for all $[x,y]\in A$ and
$t>s\ge 0$ we have
\begin{equation}\label{benilan}
\frac{1}{2}\left[\|u(t)-x\|^2-\|u(s)-x\|^2\right]\le \int_s^t\langle y,x-u(\tau)\rangle\ d\tau.
\end{equation}
Since $u$ is absolutely continuous, \eqref{benilan} implies $\dot u(t)\in -Au(t)$ almost everywhere on $[0,T]$.\\
 Monotonicity of $A$ implies   that for any proximal sequence $\{x_k\}$: $\langle x_{k-1}-x_k-\lambda_k y,x_k-x\rangle\ge
0$. But $\|x_k-x\|^2-\|x_{k-1}-x\|^2\le 2\langle x_{k-1}-x_k,x-x_k\rangle$ and so
$$\|x_k-x\|^2-\|x_{k-1}-x\|^2\le 2\lambda_k\langle y,x-x_k\rangle.$$
Summing up for $k=m+1,\dots n$ we obtain
$$\|x_n-x\|^2-\|x_0-x\|^2\le 2\sum_{k=1}^n\lambda_k\langle y,x-x_k\rangle.$$
Setting $x_0= u(s)$ and passing to the limit appropriately we finally get \eqref{benilan}. Notice that $u(t)\in D(A)$
by maximality.\bx

 A consequence of Proposition \ref{Kobayashi} and Theorem \ref{thm:proxdiff} is the following

\begin{cor}\label{cor:prox_continuous} The following statements hold:
\begin{itemize}
    \item [i)] For each $z\in D(A)$ we have
$$\|x_n-u(t)\|\le \|x_0-z\|+\|u(0)-z\|+\|A^0z\|\sqrt{(\sigma_n-t)^2+\tau_n}.$$
    \item [ii)] For trajectories $u$ and $v$ we get
$$\|v(s)-u(t)\|\le \|v(0)-z\|+\|u(0)-z\|+\|A^0z\|\ |s-t|.$$
    \item [iii)] The unique function $u$ satisfying \eqref{continuous} is Lipschitz-continuous with
$$\|u(s)-u(t)\|\le \|A^0x\|\ |s-t|.$$
    \item [iv)] $\dot u\in L^\infty(0,\infty;H)$ with $\|\dot{u}(t)\|\le\|A^0x\|$ almost
everywhere.
\end{itemize}
\end{cor}

Proposition \ref{Kobayashi} was used to construct a continuous trajectory by considering finer
and finer discretizations on a compact interval. By controlling the distance between two discrete schemes it is possible to
obtain bounds for the distance between a limit  trajectory and a discrete scheme. As a consequence, one can estimate the distance between two trajectories as well.

\subsection{Euler sequences}

Assume $A$ maps $D(A)$ into itself. (Notice that this is a strong assumption, so the range of applications of this
discretization method is limited compared to the proximal sequences). Let $\{\lambda_n\}$ be a sequence of positive
numbers or {\em stepsizes}. Define an {\it Euler  sequence} $\{z_n\}$ recursively by
\begin{equation}\label{euler}
\left\{\begin{array}{rcl}
\displaystyle \frac{z_n-z_{n-1}}{\lambda_{n-1}}  & \in  & - Az_{n-1}\qquad\hbox{for all $n\ge 1$}\\
z_0 & \in & D(A)
       \end{array}\right.
\end{equation}

\noi A remarkable feature of this scheme is that the terms of the sequence can be computed explicitly (forward scheme).\\

Observe that if $A=I-T$ with $T$ nonexpansive and $\lambda_n\equiv 1$ then $z_n=T^nz_0$. This particular case has been
studied extensively by several authors in the search for fixed points of $T$. Some of their results will be presented
in the forthcoming sections.\\
 Note also that in this framework a Kobayashi-type inequality holds too, namely
\begin{equation}\label{Eq:Kobayashi_euler}
\|z_k-\widehat z_l\|\le
    \|z_0-u\|+\|\widehat z_0-u\|+
    \|u-T(u)\|\sqrt{(\sigma_k-\widehat\sigma_l)^2+\tau_k+\widehat\tau_l},
\end{equation}
where $u$ is any point in $H$. This fact was recently pointed out by \cite[Vigeral]{Vig}.\\

Let us define the velocity at stage $n$ as
$w_n=\displaystyle\frac{ z_{n+1}-z_n}{\lambda_n}\in -Az_n.$

\begin{lem}\label{monEul}
If $y\in\S$ then $\|z_{n+1}-y\|^2\le \|z_n-y\|^2 + \lambda_n^2\|w_n\|^2$.
\end{lem}

\dem For any $y\in H$ one has
\begin{equation}\label{Eq:Euler_1}
\|z_{n+1}-y\|^2=\|z_n-y\|^2+2\lambda_n\langle w_n,z_n-y\rangle+\lambda_n^2\|w_n\|^2.
\end{equation}
The desired inequality follows from monotonicity if $0\in Ay$.\bx

Observe the similarity and the difference with \eqref{continuous} and \eqref{prox}. The dissipativity condition in Lemma \ref{monEul} is much weaker than the corresponding ones in Lemmas \ref{lem:basic_continuous_1} and \ref{lem:basic_prox_2}.\\

An immediate consequence is the following:

\begin{cor}\label{cor:basic_euler}
Assume $\sum \| z_{n+1} - z_n \|^2 <\infty$. For each $y \in \S$ the sequence $\|z_n-y\|$ is convergent.
\end{cor}

\dem It suffices to observe from Lemma \ref{monEul}  that the sequence $ \| z_n -y\| ^2 + \sum _{m=n}^{+\infty} \|z_{m+1} - z_m \|^2$ is
decreasing. \bx

{\it Comments}\\
The hypothesis in the previous result holds if $\{\lambda_n\}\in\ltwo$ and $\{w_n\}$ bounded.\theend

Notice the similarity with Corollaries \ref{cor:basic_continuous_2} and \ref{lem:basic_prox_3}.\\

The main drawback of Euler sequences is that they can be quite unstable. Most convergence results need regularity
assumptions such as $\{\lambda_n\}\in\ltwo$ and the boundedness of the sequence $\{w_n\}$, or at least that $\sum \| z_{n+1} - z_n \|^2 <\infty$.\\

An important result involving an operator $A$ of the form $I - T$ is the following, see \cite[Br\'ezis]{Bre}:

\begin{prop}[Chernoff's estimate]\label{chernoff}
Let $T$ be  non-expansive  from $H$ to itself and $\lambda >0$. If $v$ satifies
$$ \dot v(t) =-\frac{1}{\lambda}(I-T) v(t)$$
with $v(0)=v_0$ then
$$  \| v(t) - T^n v_0\| \leq \|\dot v(0)\| \sqrt{\lambda t + ( n\lambda -t) ^2}.$$
\end{prop}

\dem  It is enough to consider the case $\lambda =1$. \\Define $\phi_n(t)=\|v(t)-T^nv_0\|$ and
$\gamma_n(t)=\|\dot v(0)\|\sqrt{ t+[n-t]^2}$. We shall prove inductively that $\phi_n(t)\leq\gamma_n(t)$. For $n=0$ simply
observe that
$$\|v(t)-v_0\|\leq\int_0^t\|\dot v(s)\|\,ds\leq \|\dot v(0)\|t\leq\gamma_0(t)$$
using point 4 in Theorem \ref{exuni}.

\noi Now let us assume $\phi_{n-1}\leq\gamma_{n-1}$ and prove $\phi_n\leq\gamma_n$. Multiplying $\dot v(t)+ v(t)= Tv(t)$ by
$e^{t}$ and integrating we obtain $v(t)=v_0e^{-t}+\int_0^t e^{(s-t)}Tv(s)\,ds$ so that
\begin{eqnarray*}
\phi_n(t) & = & \left\|e^{-t}(v_0-T^nv_0)+
\int_0^t e^{(s-t)}[Tv(s)-T^nv_0]\,ds\right\|\\
& \le & e^{-t}\|v_0-T^nv_0\|+ \int_0^t e^{(s-t)}\phi_{n-1}(s)\,ds.
\end{eqnarray*}

\noi Noting that $\|v_0-T^nv_0\|\le\sum_{i=1}^n\|T^{i-1}v_0-T^iv_0\| \le n\|v_0-Tv_0\|=n \|\dot v(0)\|$ and using the
induction hypothesis we deduce
$$\phi_n(t)\leq e^{-t}\left[n\|\dot v(0)\|+
\int_0^te^{s/\lambda}\gamma_{n-1}(s)\,ds \right].$$

\noi Hence it suffices to establish the inequality
$$n+\int_0^te^{s}
\sqrt{s+[(n-1)-s]^2}\,ds \le e^{t}\sqrt{ t+[n-t]^2}.$$ Since this holds trivially for $t=0$, it suffices to prove the
inequality for the derivatives
$$e^{t/\lambda}\sqrt{ t+[(n-1)-t]^2}
\le e^{t}\left [\sqrt{t+[n-t]^2} +\frac{1-2[n-t]}{2\sqrt{ t+[n-t]^2}}\right].$$

\noi This easily verified by squaring both sides.\bx

In particular if $T$ is the resolvent  $J^A_{\lambda}$, $v$ is $u_\lambda$ and using (\ref{eqBre}),  we  deduce that
\begin{equation}\label{approxE}
\| (I + \lambda A )^{-n} x -  u(t)\| \leq  \| A^0 (u_0 ) \| \left(2\sqrt { \lambda t} + \sqrt{\lambda t + ( n\lambda
-t) ^2}\right)
\end{equation}
hence taking $\lambda = t/n$ we obtain an exponential approximation
\begin{equation}\label{approxEE}
\left\| \left(I + \frac{t}{n} A \right)^{-n} \! \!x -  u(t)\right\| \leq  3 \frac { \| A^0 (u_0 ) \|  t}{\sqrt {n}}.
\end{equation}

\subsection{ Discrete to continuous }

\noi Given a sequence $\{x_n\}$ in $X$ along with a strictly  increasing sequence $\{\sigma_n\}$
of positive numbers with $\sigma_0=0$ and $\sigma_n\to\infty$ as $n\to\infty$, one can construct a ``continuous-time"
trajectory $x$ by interpolation: for $t\in[\sigma_n,\sigma_{n+1}]$, take $x(t)$ anywhere on the  segment
$[x_n, x_{n+1}]$. It is easy to see that any trajectory defined this way converges to some $\bar x$  if, and only if,
the sequence $\{x_n\}$ converges to $\bar x$.

\noi Observe that if the interpolation is chosen to be piecewise constant in each subinterval
$[\sigma_n,\sigma_{n+1})$, then
$$\frac{1}{t}\int_0^tx(\xi)\ d\xi=\frac{1}{\sigma_n}\sum\limits_{k=1}^n\lambda_kx_k,$$
where $\lambda_k=\sigma_k-\sigma_{k-1}$. The sum on the right-hand side of the previous equality represents
an average of the points $\{x_n\}$ that is {\em weighted} by the sequence $\{\lambda_n\}$ and will be denoted by $\bar x_n$.\\

{\bf From now on we will consider only  proximal  or Euler sequences with stepsizes $\{\lambda_n \} \notin \luno$.}\\

The next sections are devoted to the asymptotic analysis. We start by considering the sequences of values in the case $A=\partial f$ in Section 3. The rest deals with the behavior of trajectories and sequences themselves. Section 4 presents general tools related to weak convergence and properties of weak limit points. These last properties are easier to satisfy for the averages and are studied in Section 5. In Section 6 we present weak convergence, in particular in the framework of demipositive operators.
Section 7 introduces different geometrical conditions that are sufficient for strong convergence. Section 8 is devoted to almost orbits and describes equivalence classes that allow to recover previous results with a new perspective and extend to non autonomous processes.


\section{Convex optimization and convergence of the values}

This section is devoted to the case  $A=\partial f$ where we evaluate $f$ on trajectories.

\subsection{Continuous dynamics}

When $A=\partial f$ with $f\in\Gamma_0(H)$ the differential inclusion  (\ref{continuous})  is a generalization of the
gradient method,  for nondifferentiable functions. In what follows let $u:[0,\infty)\to H$ be the solution of the
differential inclusion
\begin{equation}
\dot u(t)\in -\partial f(u(t)),
\end{equation}
whose existence is given in Theorem \ref{exuni}. Let $$f^* =\inf_{x\in H}f(x) \in \R \cup \{ - \infty \}.$$

The following result and its proof are essentially from \cite[Br\'ezis]{Bre} (see \cite[G\"uler]{Gul2}).\\

\begin{prop} The function $t\mapsto f(u(t))$ is decreasing and
$\lim\limits_{t\to\infty} f(u(t))= f^*$.
\end{prop}

\dem The subdifferential inequality is
$$
f(u(t))-f(u(s))\le-\langle\dot u(t),u(t)-u(s)\rangle.
$$
Thus
$$
\limsup_{s\to t^-}\frac{f(u(t))-f(u(s))}{t-s}\le-\|\dot u(t)\|^2
$$
and so the function $t\mapsto f(u(t))$ is decreasing. For each $z\in H$ and $s\in[0,t]$ the subdifferential inequality
then gives
$$
f(z)\ge f(u(s))+\langle  \dot u(s),u(s)-z\rangle\ge f(u(t))+\frac{1}{2}\frac{d}{ds}\|u(s)-z\|^2.
$$
Integrating on $[0,t]$ we obtain that
$$
tf(z)\ge tf(u(t))+\frac{1}{2}\|u(t)-z\|^2-\frac{1}{2}\|u(0)-z\|^2
$$

and so
\begin{equation}\label{ineq_f(u(t))}
f(u(t))+\frac{\|u(t)-z\|^2}{2t}\le f(z)+\frac{\|u(0)-z\|^2}{2t}
\end{equation}
for every $z\in H$.\bx

{\it Comments}\\
Inequality \eqref{ineq_f(u(t))} shows that if $\S\neq \emptyset$ then $f(u(t))$ converges to $f^*$ at a rate of
$O(1/t)$. However, if the trajectory $u(t)$ is known to have a strong limit, then the rate drops to $o(1/t)$ (see
\cite[G\"uler]{Gul2}).

 \theend

\subsection{Proximal sequences}
Let $\{x_n\}$ be a proximal sequence associated to $A = \partial f$. The following result is due to \cite[G\"uler]{Gul}:

\begin{prop}
The sequence $f(x_n)$ is decreasing and $\limn f(x_n)=f^*$.
\end{prop}

\dem The subdifferential inequality implies $f(x_{n-1}) - f(x_n) \geq \lambda_n\|y_n\|^2$ so that $f(x_n)$ is
decreasing. Convergence of $f(x_n)$ to $f^*$ follows from Lemma \ref{lem:Gulsubdif} below since $\sigma_n\to\infty$.\bx

\begin{lem}\label{lem:Gulsubdif} Let $u\in\dom f$, then
$$f(x_n) - f(u) \leq \frac{\| u-x_0\|^2}{2\sigma_n} - \frac{\|u - x_n\|^2}{2\sigma_n} - \frac{\sigma_n}{2} \|y_n\|^2.$$
\end{lem}

\dem The subdifferential inequality  is
$$f(u) - f(x_n) \geq \langle u-x_n, -y_n\rangle = \frac { \langle u- x_n, x_{n-1} - x_n \rangle} { \lambda_n}
$$
for all $u$ in the domain of $f$. Thus

$$2{ \lambda_n} (f(u) - f(x_n)) \geq  \| u- x_n\|^2 +  \lambda_n^2\|y_n\|^2 - \|u -x_{n-1}\|^2.$$
Summing  up from $1$ to $n$ leads to
\begin{equation}\label{gulvalue}
2\sigma_n f(u) - 2\sum_{k=1}^n{ \lambda_k}f(x_k) \geq  \| u- x_n\|^2 + \sum_{k=1}^n \lambda_k^2\|y_k\|^2 - \| u-
x_0\|^2.
\end{equation}

\noi On the other hand the subdifferential inequality implies $f(x_{n-1}) - f(x_n) \geq \lambda_n\|y_n\|^2$. Multiplying by $\sigma_{n-1}$ and
rearranging we get
$$ \sigma_{n-1}f(x_{n-1}) - \sigma_n f(x_n) + \lambda_n f(x_n) \geq \lambda_n\sigma_{n-1} \|y_n \|^2,$$
from which we derive
$$- \sigma_n f(x_n) + \sum_{k=1}^n \lambda_k  f(x_k) \geq \sum_{k=1}^n\lambda_k\sigma_{k-1} \|y_k \|^2$$
by summation. Adding twice this inequality to \eqref{gulvalue} we obtain
$$2\sigma_n ( f(u) - f(x_n)) \ge
\|u- x_n\|^2 - \|u- x_0\|^2 +  \sum_{k=1}^n\lambda_k^2\|y_k\|^2 + 2  \sum_{k=1}^n\lambda_k\sigma_{k-1}\|y_k\|^2.$$
Recall from Lemma \ref{lem:basic_prox_1} that $\|y_n\|$ is decreasing. We get
$$\|y_n\|^2  \sigma_n^2=\|y_n\|^2\sum_{k=1}^n(\lambda_k^2+ 2\lambda_k\sigma_{k-1})\le\sum_{k=1}^n(\lambda_k^2+ 2\lambda_k\sigma_{k-1})\|y_k\|^2$$
and the result follows at once by rearranging the terms. \bx

{\it Comments}\\
\noi If $\S\neq\emptyset$, Lemma \ref{lem:Gulsubdif} gives
\begin{equation}\label{proxspeed}
\|y_n\|\le \frac{d(x_0,\S)}{\sigma_n}.
\end{equation}
A similar estimation had been proved in \cite[Br\'ezis and Lions]{BrL} but the right-hand side is $\sqrt{2}$ times
larger.\\

The fact that $f(x_n)\to f^*$ had first been proved in \cite[Martinet]{Mar} when $f$ is
coercive and $\lambda_n\equiv\lambda$.\\

By Lemma \ref{lem:Gulsubdif}, if $\S\neq\emptyset$ the rate of convergence can be estimated at $O(1/\sigma_n)$.
Moreover, \eqref{proxspeed} and the subdifferential inequality together give
$$f(x_n)-f^* \le \langle x^*-x_n, -y_n\rangle\le \|x^*-x_n\|\ \|y_n\|\le \frac{d(x_0,\S)\|x^*-x_n\|}{\sigma_n}$$ for
all $x^*\in\S$. Therefore, if the sequence $\{x_n\}$ is known to converge strongly, then $|f(x_n)-f^*|=o(1/\sigma_n)$.
This was proved in \cite[G\"uler]{Gul} using a clever but unnecessarily sophisticated argument instead of inequality
\eqref{proxspeed}. \theend

\subsection{Euler sequences}

In this case the sequence $f(z_n)$ need not be decreasing. However, we have the following:

\begin{lem}\label{lem:Euler_liminf}
If either $\sum\|z_{n+1}-z_{n}\|^2<\infty$ or $\limn\lambda_n\|w_n\|^2=0$, then $\liminf\limits_{n\to\infty}f(z_n)=f^*$.
\end{lem}

\dem Since $-w_n\in\partial f(z_n)$, the subdifferential inequality and \eqref{Eq:Euler_1} together imply
\begin{equation}\label{Eq:Euler_liminf}
\|z_{n+1}-y\|^2\le\|z_n-y\|^2-2\lambda_n(f(y)-f(z_n))+\lambda_n^2\|w_n\|^2
\end{equation}
for each $y\in H$. If $\sum\|z_{n+1}-z_{n}\|^2<\infty$ then
$$\sum\lambda_n(f(z_n)-f(y))<\infty$$
(possibly $-\infty$). Since $\{\lambda_n\}\notin\luno$ one must have $\liminf\limits_{n\to\infty}f(z_n)\le f(y)$ for each $y\in H$.\\

On the other hand, inequality \eqref{Eq:Euler_liminf} can be rewritten as
$$\lambda_n\left[2\frac{}{}\!(f(z_n)-f(y))-\lambda_n\|w_n\|^2\right]\le \|z_n-y\|^2-\|z_{n+1}-y\|^2$$
so that
$$\sum\lambda_n(f(z_n)-f(y)-\lambda_n\|w_n\|^2)<\infty$$
and $\liminf\limits_{n\to\infty}f(z_n)\le f(y)$ for each $y\in H$.\bx

A complementary result is the following from \cite[Shor]{Sho}:

\begin{prop}
Let $\hbox{dim}(H)<\infty$ and assume $\S$ is nonempty and compact. If $\limn\lambda_n=0$ and the sequence $w_n$ is bounded
then $\limn f(z_n)=f^*$.
\end{prop}

\dem By continuity, it suffices to prove that $\hbox{dist}(z_n,\S)=\inf_{y\in\S}\|z_n-y\|$ tends to $0$ as $n\to\infty$. For $\gamma>f^*$ define $L_\gamma=\{x:f(x)=\gamma\}$. Take $\eps>0$ and define
$$\delta(\eps)=\min_{y\in\S}\hbox{dist}(y,L_{f^*+\eps})\qquad\hbox{and}\qquad d(\eps)=\max_{y\in\S}\hbox{dist}(y,L_{f^*+\eps}).$$
Observe that $0<\delta(\eps)\le d(\eps)\to 0$ as $\eps\to 0$. By hypothesis and Lemma \ref{lem:Euler_liminf} there is $N\in\N$ such that $f(z_N)\le f^*+\eps$ and $\lambda_n\|w_n\|\le\delta(\eps)$ for all $n\ge N$.
We shall prove that $\hbox{dist}(z_n,\S)\le 2d(\eps)$ for all $n\ge N$. Since $\eps>0$ is arbitrary this shows that $\limn\hbox{dist}(z_n,\S)=0$.\\

Indeed, if $f(z_n)\le f^*+\eps$ (this holds for $n=N$) then $\hbox{dist}(z_n,\S)\le d(\eps)$ and $\hbox{dist}(z_{n+1},\S)\le d(\eps)+\delta(\eps)\le 2d(\eps)$. On the other hand, if $f(z_n)>f^*+\eps$ then $\hbox{dist}(z_{n+1},\S)\le\hbox{dist}(z_n,\S)$. To see this, notice that if $y\in\S$ then $\langle \frac{w_n}{\|w_n\|},y-z_n\rangle$ is the distance from $y$ to the hyperplane $\Pi_n=\{x:\langle w_n,z_n-x\rangle\}$, so that
$$\langle w_n,y-z_n\rangle  \ge  \|w_n\|\hbox{dist}(\S,\Pi_n) \ge  \|w_n\|\hbox{dist}(\S,L_{f(z_n)}) \ge  \|w_n\|\delta(\eps),$$
where the second inequality follows from convexity and the last one is true whenever $f(z_n)>f^*+\eps$. Using \eqref{Eq:Euler_1} and recalling that $\lambda_n\|w_n\|\le\delta(\eps)$ we deduce that
$$\hbox{dist}(z_{n+1},\S)^2\le\hbox{dist}(z_n,\S)^2-\lambda_n\|w_n\|\delta(\eps),$$
proving that $\hbox{dist}(z_{n+1},\S)\le\hbox{dist}(z_n,\S)$.\bx

Observe that this result does not require the stabilizing summability condition but it is necessary to make a very strong
assumption on the set $\S$.\\

\section{General tools for weak convergence}

We denote by $\Omega [u(t)]$ (resp. $\Omega [x_n]$) the set of weak cluster points of a trajectory $u(t)$ as $t\to\infty$
(resp. of a sequence $\{x_n\}$ as $n\to\infty$). \\

Given a  trajectory $u(t)$ we define
$$ \bar u (t)  = \frac{1}{t}\int_0^tu(\xi)\ d\xi$$
Similarly, given  a sequence $\{x_n\}$ in $H$ along with stepsizes $\{\lambda_n\}\notin\luno$, we introduce
 $$\bar x_n = \frac{1}{\sigma_n}\sum_{k=1}^n\lambda_kx_k.$$

\subsection{Existence of the limit}

Most of the results on weak convergence that exist in the literature rely on the combination of two
types of properties involving a subset $F \subset H$:\\

The first one is a kind of ``Lyapounov condition" on the sequence or the trajectory like
\begin{itemize}
    \item [{\bf (a1)}] $\|x_n -u \|$ converges to some $\ell (u)$ for each $u \in F$, or
    \item [{\bf (a2)}] $P_F (x_n)$ converges strongly (in all that follows $F$ will be closed and convex).
\end{itemize}
These properties imply that the sequence is somehow ``anchored" to the set $F$.\\

The second one is a global one, concerning the set of weak cluster points of the sequence or trajectory:
\begin{itemize}
    \item [{\bf (b)}] $\Omega [x_n] \subset F$.
\end{itemize}
However, it is sometimes available only for the averages:
\begin{itemize}
    \item [{\bf (b')}] $\Omega [\bar x_n] \subset F$.
\end{itemize}

\noi The following result is a very useful tool for proving weak convergence of a sequence
on the basis of {\bf (a1)} and {\bf (b)} above. It is known, especially in Hilbert spaces, as {\em Opial's lemma} \cite{Opi}.


\begin{lem}[Opial's Lemma]\label{Op}
Let $\{x_n\}$ be a sequence in $H$ and let $F\subset X$. Assume
\begin{enumerate}
   \item $\| x_n - u\|$ has a limit as $n\to\infty$ for each
$u\in F$; and
  \item $\Omega[x_n] \subset F$.
\end{enumerate}
 Then $x_n$ converges weakly to some $x^*\in F$.
\end{lem}
\dem Since $\{x_n\}$ is bounded it suffices to prove that it has only one weak cluster point. Let $x, y \in  \Omega[x_n] \subset F$
so that $\| x_n - x\|$ converges to $\ell (x)$ and similarly for $y$. From
$$\|x_n-y\|^2=\|x_n-x\|^2+\|x-y\|^2+2\langle x_n-x,x-y\rangle$$
one deduces by choosing appropriate subsequences
$$
\ell(y) = \ell (x)  +  \| x-y \|^2    \qquad (x_{\phi(n)}\wto x)
$$ and
$$
\ell(y) = \ell (x)  -  \| x-y \|^2    \qquad (  x_{\psi(n)}\wto y)
$$
hence $x=y$. \bx

{\it Comments}\\
A Banach space $X$ satisfies {\em Opial's condition} if it is reflexive and
\begin{equation}\label{Eq:Opial}
\limsup\limits_{n\to\infty}\|x_n-x\|<\limsup\limits_{n\to\infty}\|x_n-y\|\quad\hbox{whenever}\quad x_n\wto x\neq
y.\end{equation}holds. Any uniformly convex Banach space having a weakly continuous duality mapping (in particular, any
Hilbert space) satisfies Opial's condition (see \cite[Opial]{Opi}). Opial's Lemma holds in any Banach space satisfying
Opial's condition. \theend

Following \cite[Passty]{Pas1}, one obtains a more general result:

\begin{lem}\label{A1}
Let $\{x_n\}$ be a sequence in $H$ with stepsizes $\{\lambda_n\}$  and let $F\subset X$. Assume {\bf(a1)} : the
sequence $\| x_n - u\|$ has a limit as $n\to\infty$ for each $u\in F$. Then the sets $\Omega[x_n] \cap F$ and
$\Omega[\bar x_n] \cap F$ each contains at most one point. In particular if $\Omega[x_n]\subset F$ (resp. $\Omega[\bar
x_n]$), then $x_n$ (resp. $\bar x_n$) converges weakly as $n\to\infty$. A similar result holds for trajectories.
\end{lem}

\dem Write
$$\|x_n-y\|^2=\|x_n-x\|^2+\|x-y\|^2+2\langle x_n-x,x-y\rangle$$
So that $\langle x_n,x-y\rangle$ converges to some $m(x,y)$ for any $x, y \in F$. If $u$ and $v$ belong to $ \Omega[
x_n] \cap F$  one obtains  $\langle u, u  - v \rangle =   \langle v, u  - v \rangle$  hence $u= v$. Similarly $\langle
\bar x_n,x-y\rangle$ converges to  $m(x,y)$. Thus   both $ \Omega[ x_n] \cap F$ and  $ \Omega[\bar x_n] \cap F$ contain
at most one point.\bx

An alternative proof using {\bf (a2)} and either {\bf (b)} or {\bf (b')} is as follows:

\begin{lem}\label{A2}
Let $\{x_n\}$ be a  bounded sequence in $H$  with stepsizes $\{\lambda_n\}$ and let $F\subset X$ be closed and convex.
Assume  {\bf(a2)}:  $P_F x_n \to \zeta$  as $n\to\infty$. Then
$$\Omega[x_n]\cap F = \Omega[\bar x_n]\cap F = \{\zeta\}.$$
In particular, if $ \Omega[x_n]\subset F$ (resp. $ \Omega[\bar x_n]$), then $x_n$ (resp. $\bar x_n$) converges weakly
to $\zeta$. A similar result is true for trajectories.
\end{lem}

\dem By definition of the projection, for each $u\in F$ one has
$$\langle x_n - P_Fx_n, u - P_F x_n \rangle  \leq 0.$$
Since $x_n$ is bounded we deduce that
$$\langle x_n - \zeta, u - \zeta \rangle  \leq \rho_n$$
with $\limn\rho_n=0$. This implies $ \Omega [ x_n] \cap F = \{\zeta \}$ (if $v\in \Omega [ x_n] \cap F$, take $u = v$).
Similarly
$$\langle \bar x_n - \zeta, u - \zeta \rangle  \leq \bar\rho_n,$$
which gives $ \Omega [\bar  x_n] \cap F = \{\zeta \}$.\bx

In our case the set $F$ will always be $\S$, which is closed and convex.\\

\subsection{Characterization of the limit: the asymptotic center}

We show here that moreover the weak limit can be characterized.\\

Given a bounded  sequence $\{x_n\}$ let
$$G(y) = \limsup_{n\to\infty} \|x_n - y \|^2$$
(for a trajectory $u(t)$ define $G(y)=\limsup\limits_{t\to\infty}\|u(t) - y \|^2$). The function $G(y)$ is continuous,
strictly convex and coercive. Its unique minimizer is called the {\em asymptotic center} (see \cite{Ede}) of the
sequence (resp. trajectory) and is denoted by $AC\{ x_n\}$
(resp. $AC\{ u(t)\}$).\\

Observe that, by virtue of Opial's condition \eqref{Eq:Opial}, if $x_n\wto x$ then $x=AC\{ x_n\}$.\\

The weak limit of the average is still the asymptotic center, under some assumptions.

\begin{prop}\label{AC}
Assume {\bf (a1)}. If $\bar x_n \wto x \in F$, then $x = AC\{ x_n\}$.\\
The same property holds for trajectories.
\end{prop}

\dem For each $y \in H$ we have
$$\|x_n -x\|^2 = \| x_n - y\|^2 + 2 \langle x_n - y , y - x  \rangle + \| y - x \| ^2.$$
Hence
$$\frac{1}{\sigma_n}\sum_{m=1}^n\lambda_m\|x_m -x\|^2 =
\frac{1}{\sigma_n}\sum_{m=1}^n\lambda_m\|x_m - y\|^2 + 2\langle \bar x_n-y,y-x\rangle + \|y-x\| ^2.$$

If $x_n \wto x$ and $x\in F$ then $\ell(x)=\limn\|x_n-x\|$ exists. Therefore,
$$G(x)=\ell(x)^2 \leq \limsup_{n\to\infty}\left[ \frac{1}{\sigma_n}\sum_{m=1}^n\lambda_m\|x_m - y\|^2\right] - \| x-y \|^2 \le G(y) - \| x-y \|^2$$
for each $y\in H$ so that $x = AC\{ x_n\}$.\bx

\subsection{Characterization of the weak convergence}

In this section we use the fact that the trajectories or sequences are generated through a maximal monotone operator.\\

Let us consider first the case $A = I-T$,  where $T$ is non expansive. The following result is in \cite[Pazy]{Paz3}:

\begin{prop}\label{PP1}
The sequence $T^nx$ converges weakly if, and only if, $\S \neq \emptyset$ and $\Omega [T^n x] \subset \S$.
\end{prop}

\dem Assume $\S \not= \emptyset$.  Given $u \in \S$, the sequence
$\|T^n x - u \|$ is decreasing and so $T^n x $ is bounded. By Lemma \ref{A1}, the fact that $\Omega [T^n x] \subset \S$
implies that $T^nx$ converges weakly. Conversely, since the sequence $\{T^n x\}$ is bounded, the argument in the  proof  of  Theorem \ref{thm:Bai_T_average}
shows that the weak limit of $T^nx$ must be in $\S$.\bx

An alternative proof relies on the following result, which is interest in its own right:

\begin{lem}
Assume the sequence $U_n x  = \frac{1}{n}(z + Tz + ... + T^{n-1}z)$ is bounded. Then $\emptyset \neq \Omega [U_n x ]
\subset \S$.
\end{lem}
\dem For any $y\in H$ one has
\begin{eqnarray*}
0  &  \leq  &  \| T^k x - y \|^2 - \| T ^{k+1}x - T y \|^2\\
&  =  &  \| T^k x - Ty \| ^2- \| T ^{k+1}x - T y \|^2  + \| Ty - y \| ^2 + 2 \langle T^k x - T y , Ty - y \rangle.
\end{eqnarray*}
By taking the average we obtain
$$0 \leq \frac{1}{n} \| x - Ty\|^2  + \| Ty - y \| ^2  +2  \langle U_n x - T y , Ty - y \rangle.$$
Therefore, if $p\in \Omega [U_n x ]$, we can let $n\to\infty$ to deduce that
$$ 0 \leq  \| Ty - y \| ^2  +2  \langle p - T y , Ty - y \rangle.$$
In particular, if $y = p$ we conclude that $ \|Tp -p \| ^2 \leq 0 $ and so $p\in \S$. \bx

Assuming that $\S$ is nonempty we can give a direct proof:

\begin{lem}
Assume $\S \neq \emptyset$. Then $T^nx \wto p$ implies $ p \in \S$.
\end{lem}

\dem For any $y\in H$ and $u\in \S$
\begin{eqnarray*}
0  &  \leq  &  \| T^k x - y \|^2 - \| T ^{k+1}x - T y \|^2\\
&  =  &  \| T^k x - u  \|^2 - \| T ^{k+1}x - u \|^2  + \|u - y \|^2 - \| u - Ty \|^2\\
&     &  + 2 \langle T^k x -u, u-y\rangle -2 \langle T^ {k+1} x -u, u- Ty\rangle.
\end{eqnarray*}
Take $y = p$ and let $k\to\infty$. Since $\lim\limits_{k\to\infty}\|T^kx-p\|$ exists we get
$$
0 \leq \|u - p \|^2 + 2 \langle p -u, u-p\rangle -  \| u - Tp \|^2  - 2 \langle p -u, u- Tp\rangle.
$$
which is precisely $\|Tp -p \| ^2 \leq 0$ and implies $p\in\S$.\bx

Following \cite[Pazy]{Paz4}, one obtains the continuous counterpart of Proposition \ref{PP1}:

\begin{prop}\label{PP2}
The trajectory $S_tx$ converges weakly if, and only if, $\S \neq\emptyset$ and $\Omega [S_t x] \subset \S$.
\end{prop}

\dem Assume $\S \not = \emptyset$. By Corollary  \ref{cor:basic_continuous_2}  and Lemma \ref{A1}, $\Omega [S_t x]
\subset \S$ implies $S_t x$ converges weakly.\\ It remains to prove that if  $S_tx \wto y$ then $y \in \S$. To see
this, take any $[u, w ] \in A$. We have
\begin{eqnarray*}
\| S_t x - u \|^2  - \| x - u \|^2  &  \leq  &  2 \int _0 ^t \langle w, u - S_s x\rangle ds\\
&  =  &  2t \langle w, u - y\rangle  + 2 \int _0 ^t \langle w,y - S_s x\rangle\ ds.
\end{eqnarray*}
It suffices to divide by $t$ and let $t \to \infty$ to obtain
$$0 \leq \langle w, u - y\rangle$$
so that $y\in \S$ by maximality. \bx

Note that the proof uses the generator $A$ (compare to the  proof of the previous Proposition \ref{PP1}).\\

A last result, due to \cite [Bruck]{Bru5}, shows that if $\S\neq\emptyset$, then weak convergence is equivalent to weak
asymptotic regularity. We follow \cite[Pazy]{Paz5}.

\begin{prop}\label{PP3}
Assume $\S\neq\emptyset$. The trajectory $S_tx$ converges weakly if, and only if,
$$S_{t+h}x  - S_t x  \wto 0 \quad  \mbox{as \ }  t \to \infty$$
for each $h\ge 0$. A similar result holds for the sequence $T^n x$.
\end{prop}

\dem For $u \in \S$ and $t>s$ we have
\begin{center}$2\langle S_{s+h}x - u, S_s x -u \rangle -2\langle S_{t+h}x - u, S_t x -u \rangle\leq    \| S_{s+h} x - u \|^2  -
\| S_{t+h} x - u \|^2+  \| u - S_s x \|^2 - \|u  - S_t x \|^2.$\end{center}

Let $w \in \Omega[S_t x ]$ and $h_k \to \infty$ with $S_{t+ h_k} \wto w$. Then  $S_{s+ h_k} \wto w$ as well by weak
asymptotic regularity. Thus we obtain
$$
2 \langle w - u, S_s x -S_t x  \rangle
\leq
  \| u - S_s x \|^2 - \|u  - S_t x \|^2.
$$
so that  by {\bf(a1)}, $\langle w - u, S_t x  \rangle $ has a limit $L(w)$. In particular  $w' \in  \Omega[S_t x ]$ implies
$\langle w - u, w'  \rangle = L(w) $ so that $\langle w - u, w' - w \rangle = 0 $.
Hence by symmetry  $\langle w' - u, w - w' \rangle = 0 $, thus  $w=w'$ and $ \Omega[S_t x ]$ is reduced to one point.
 \bx

\section{Weak convergence in average}

A trajectory $u(t)$ {\em converges in average} if
$$\bar u(t)=\frac{1}{t}\int_0^tu(\xi)\ d\xi\qquad\hbox{converges as}\quad t\to\infty.$$
Similarly, consider a sequence $\{x_n\}$ in $H$ along with stepsizes $\{\lambda_n\}$,
 then $\{x_n\}$ {\em converges in average} if
 $$\bar x_n = \frac{1}{\sigma_n}\sum_{k=1}^n\lambda_kx_k \qquad\hbox{converges as}\quad n\to\infty.$$

\subsection{Continuous dynamics}

Consider $x\in\overline{D(A)}$. In order to use the semigroup notation, let us introduce
$$
\sigma_tx = \frac{1}{t} \int _0^t S_s x\ ds.\footnote{More generally $\sigma _n x = \int _0^\infty S_sx\ a_n(s)\ ds$
where $a_n$ is the density of a positive probability measure on $\R^+$, which is assumed to be of bounded variation
with $\int _0^\infty |da_n| \rightarrow 0$.}$$

\noi In order to prove that $\sigma_tx$ converges weakly as $t\to\infty$ we follow the ideas in \cite[Baillon and
Br\'ezis]{BaB}.  We first prove that the projection $P_{\S}S_tx$ converges  strongly  to some $v$  {\bf (a2)}, next
that weak cluster points of $\sigma_tx$ are in $\S$ {\bf (b')}, and finally use Lemma \ref{A2} to conclude that
$\sigma_tx$ converges weakly to $v$.

\begin{lem}\label{cvproj}
Assume $\S\neq\emptyset$. Then $ P_{\S}S_tx$ converges strongly.
\end{lem}

\dem Let  $v(t) =P_{\S}S_tx$ and  observe that the function $\psi(t)=\|v(t) - S_tx\|$ is decreasing:
$$\psi(t+h)\leq \|v(t) - S_{t+h}x \| = \|S_hv(t) - S_hS_tx \|  \leq \psi(t).$$
Therefore, it has a limit as $t\to\infty$. On the other hand, the parallelogram equality gives
$$
\|v(t+h)-v(t)\|^2+4\left\|\sfrac{v(t+h)+v(t)}{2}-S_{t+h}x\right\|^2 =
2\|v(t+h)-S_{t+h}x \|^2 + 2\|v(t)-S_{t+h}x \|^2.
$$

\noi  $\S$  convex implies $\left\|\frac{v(t+h)+v(t)}{2}-S_{t+h}x\right\|\ge\psi(t+h)$. We finally get
$$\|v(t+h) - v(t)\|^2  \leq  2 \left[\psi(t)^2 - \psi(t+h)^2\right]$$
and conclude that $v(t)$ has a strong  limit $v$ as $t\to\infty$.\bx

\begin{lem}\label{almostcluster}
 $\Omega[\sigma_tx] \subset \S$.
\end{lem}
\dem Assume $\sigma_{t_k}x\wto u$ as $k\to\infty$ and recall that  $u(t) = S_tx$. For any $v\in D(A)$ we have
$$2\int_0^{t_k} \langle u(t) - v, \dot u(t) \rangle\  dt = \|u(t_k) - v\|^2-\|x - v\|^2.$$

\noi Now take $w\in Av$, so that  $ \langle u(t) -v,  -w \rangle\geq
\langle u(t) -v, \dot u(t) \rangle$. This gives
$$2 \int_0^{t_k} \langle w, v -S_tx \rangle\ dt \geq \| S_{t_k}x - v \|^2 - \|x- v\|^2\ge - \|x- v\|^2.$$
Divide by $t_k$ and take the weak limit as $k\to\infty$. We get $\langle w, v -u \rangle \geq 0$ for any $[v, w] \in
A$, so $0 \in Au$ by maximality.\bx

{\it Comments}\\
\noi Lemma \ref{almostcluster} implies that if $\S=\emptyset$ then $\|\sigma_tx\|\to\infty$ for every
$x\in\overline{D(A)}$ as $t\to\infty$. On the other hand, if
$\S\neq\emptyset$ then every trajectory $S_tx$ is bounded, so $\sigma_tx$ is bounded for all $x\in\overline{D(A)}$.
\theend

Using Lemma \ref{A2},  Lemma \ref{cvproj} and Lemma \ref{almostcluster} we finally obtain

\begin{thm}\label{contaverage}
If $\S\neq\emptyset$,  then $\sigma_tx$ converges weakly to $v=\limty{t}P_{\S}S_tx$.
\end{thm}

As a consequence of Proposition \ref{AC} one has
\begin{prop}
If $\S\neq\emptyset$, the limit $w-\limty{t}\sigma_tx$ is the asymptotic center $AC \{ S_tx \}$.
\end{prop}

\noindent{\it Comments}\\
Weak convergence in average is still true in uniformly convex Banach space with Fr\'echet-differentiable norm (see
\cite[Reich]{Rei}) or satisfying Opial's condition (see \cite[Hirano]{Hia}).\theend

\noi If $A=\partial f$ with $f\in\Gamma_0(H)$, convergence in average guarantees the convergence of the trajectory (see
\cite[Bruck]{Bru1}):

\begin{prop}If $A=\partial f$ then
$\limty{t} \left\|u(t) - \frac {1}{t} \int_0^t u(s)\ ds\right\| = 0$.
\end{prop}

\dem Integration by parts gives $u(t)- \frac {1}{t} \int_0^t u(s)\ ds= \frac {1}{t} \int_0^t s \dot u(s)\ ds.$ $\|\dot
u(t)\|$ being decreasing  by Proposition \ref{u'dec},  one has
$$
\int_{t/2}^ts\|\dot u(s)\|^2\ ds\ge\|\dot u(t)\|^2\int_{t/2}^ts\ ds=\frac{3}{8}t^2\|\dot u(t))\|^2.
$$
 But in the case  $A=\partial f$, the function $t\mapsto t\|\dot u(t)\|^2$ is in $L^1(0,\infty)$
(see \cite[Br\'ezis]{Bre1}) which implies $\limty{t}t\|\dot u(t)\|=0$ and the result follows.\bx

It is known that both the trajectory and the average converge weakly (Theorems \ref{contaverage} and \ref{demiweak}). The preceding result implies, in particular, that the
average cannot converge strongly unless the trajectory itself does.\\

\subsection{Proximal sequences}

Consider a proximal sequence $\{x_n\}$ in $H$ along with stepsizes $\{\lambda_n\}$,
and recall that  $\bar x_n = \frac{1}{\sigma_n}\sum_{k=1}^n\lambda_kx_k$.

The next result was presented in \cite[Lions]{Lio}:

\begin{thm}\label{proxav}
Let $\S\neq\emptyset$. Then $\{x_n\}$ converges weakly in average to a point in $\S$.
\end{thm}

\dem The case $\{\lambda_n\}\notin\ltwo$ will follow from Theorem \ref{weakprox}, which states that $\{x_n\}$ converges
weakly under this condition. Therefore we assume $\{\lambda_n\}\in\ltwo$ and  check the conditions of Lemma \ref{A1}
with $F=\S$: {\bf (a1)} follows from Corollary \ref{lem:basic_prox_3}, while {\bf (b')} follows from Lemma
\ref{almostclusterprox} below. \bx

\begin{lem}\label{almostclusterprox}
Assume $\{\lambda_n\}\in\ltwo$, then $\Omega[ \overline{x}_n]\subset\S$.
\end{lem}
Take $[u,v]\in A$ and use (\ref{incnorm}) so that
$$
\|u-x_{n+1}\|^2\le\|u-x_n+\lambda_nv\|^2=\|u-x_n\|^2+\lambda^2_n\|v\|^2-2\langle v, \lambda_nx_n-\lambda_nu \rangle.
$$
Summing up for $k=1,\ 2,\ \dots\ n$ and dividing by $\sigma_n$ we obtain
$$
2\langle v,\overline{x}_n-u\rangle\le\frac{1}{\sigma_n}\|x_0-u\|^2+\frac{\tau_n}{\sigma_n}\|v\|^2.
$$
If $\overline{x}_n\wto\overline{x}$,  then
$\langle v, u - \overline{x}\rangle\ge 0$, hence $\overline{x}\in\S$ by maximality. \bx

This is the counterpart of Lemma \ref{almostcluster}.\\

The extension to the sum of two operators is in \cite[Passty]{Pas1}.\\

\subsection{Euler sequences}

For nonexpansive mappings,  weak convergence in average of the discrete iterates was established in \cite[Baillon]{Bai1}. The proof is again of the form {\bf (a2)} and {\bf (b')} but note that the property $\S \not= \emptyset$ is not assumed but  obtained during  the proof.

\begin{thm}\label{thm:Bai_T_average}
Let $T$ be a nonexpansive mapping on a bounded closed convex subset $C$ of $H$. For every $z \in C$ the sequence
$z_n=T^nz$ converges weakly in average to a fixed point of $T$, which is the strong limit of the sequence $P_{\S}T^nz$.
\end{thm}

\dem
Note that for any  $a$ and $a^i, i = 0, ..., n-1,$ in $H$, the quantity
$$
\left\| a - \frac{1}{n} \sum_{i=0}^{n-1}a^i\right\|^2 - \frac{1}{n}  \sum_{i=0}^{n-1}\| a - a^i \|^2
$$
 is independent of $a$. Hence with $U_nz = \frac{1}{n}(z + Tz + ... + T^{n-1}z)$ one has
$$
\| T U_n z - U_n z\| ^2 = \frac{1}{n}  \sum_{i=0}^{n-1} \| T U_n z - T^i z \| ^2 -  \frac{1}{n}  \sum_{i=0}^{n-1} \| U_n z - T^i z \|^2
$$
$$
\leq  \frac{1}{n} \left( \| T U_n z - z \|^2 - \| U_n z - T^{n-1} z \| ^2\right)
$$
so that
$$
\| T U_n z - U_n z\| \leq  \frac{1} {\sqrt n}  \| T U_n z - z \|.
$$
Thus  $T U_n z - U_n z \to 0$ and if $ U_n z \wto u $ then $Tu = u $ by Proposition \ref{prop:ws_sw}. It follows that $ \Omega[ U_n z] \subset  \S$, which is {\bf (b')} and $\S \not=\emptyset$.
Since, for $ u \in \S$, $ \|T^n z - u\|$ decreases,  then letting  $V_n z = P_{\S} T^n z$, $\| T^n z  - V_n z\|  $  decreases as well,  hence $V_n z$ converges to some $V$ (like in the proof of Lemma  \ref{cvproj})
 which implies that  $\Omega[ U_n z] = \{ V\}$ by Lemma \ref{A2}.
\bx

{\it Comments}\\
The conclusion of Theorem \ref{thm:Bai_T_average} holds also if $X$ is uniformly convex with
Fr\'echet-differen\-tia\-ble norm and $\lambda_n\to 1$ or if $X$ is superreflexive (\cite[Reich]{Rei}). \theend

By following an idea of Konishi (see \cite[Baillon]{BaiTT}) one can prove that the ergodic
theorem for non expensive mappings implies in fact the analogous results for the semi-group:

\begin{prop}
Theorem \ref{thm:Bai_T_average} implies Theorem \ref{contaverage}.
\end{prop}

\dem Let $0<h<t$ and $n= [ t/h]$ the integer part of $t/h$ and set $T_h = S_h$ and $U_nx = \frac {1}{n}
\sum_{m=0}^{n-1} T ^m x$. One has
$$
t \sigma_t x =  \int_0^h S_s x ds +  ... + \int_{(n-1)h}^{nh} S_sx ds + \int_{nh}^t S_s x ds
$$
and
$$
\left\|  \int_0^h S_s x  ds - hx \right\| \leq \int_0^h \|S_s x- x \| ds.
$$
Similarly
$$
\left\|  \int_{mh}^{(m+1)h} S_s x  ds - h T_h^mx \right\| \leq \int_{mh}^{(m+1)h} \|S_s x-  S_{mh}x \| ds\leq \int_0^h \|S_s x- x \| ds
$$
hence
$$
\|t \sigma_t x - nh U_nx \| \leq n \int_0^h \|S_s x- x \| ds +M h,
$$
where $\|S_s x \| \leq M$. Thus
$$
\| \sigma_t x -  U_nx \| \leq \frac{1}{h} \int_0^h \|S_s x- x \| ds + \frac{2M}{n}.
$$
But as $t \to + \infty$, $U_nx$ converges weakly to a fixed point $u_h$ of $T_h$ by  Theorem \ref{thm:Bai_T_average}. \\
Let us now  prove that $u_h$ is a Cauchy net  as $h \to 0$. Given $0<h, h' < t$, $n = [t/h], n' = [ t/h']$ one has
$$
\|U_n x - U_{n'} x \| \leq  \frac{1}{h} \int_0^h \|S_s x- x \| ds + \frac{2M}{n} +  \frac{1}{h'} \int_0^{h'} \|S_s x- x \| ds + \frac{2M}{n'}.
$$
Hence as $t \to + \infty$
$$
\|u_h - u_{h'} \| \leq  \frac{1}{h} \int_0^h \|S_s x- x \| ds + \frac{1}{h'} \int_0^{h'} \|S_s x- x \| ds,
$$
thus  $u_h$ is a Cauchy net  that  converges to some $u$, since $\| S_sx -x \| \to 0$ as $s \to 0$. But $S_{mh} u_h= u_h$, so that  given $s$ and $h = s/m$ one has $S_s u_h = u_h$.  As $m \to + \infty$  this implies $S_s u = u$, thus   $u\in \S$. Now write, given $y\in H$
$$
| \langle \sigma_t x  - u, y \rangle | \leq | \langle\sigma_t  x- U_n x , y \rangle | + |  \langle U_n x - u_h, y \rangle | + \| u_h - u\| \| y\|
$$
hence
$$
| \langle \sigma_t x  - u, y \rangle | \leq  \left(  \frac{1}{h} \int_0^h \|S_s x- x \| ds + \frac{2M}{n} \right)
\|y \| + |  \langle U_n x - u_h, y \rangle | + \| u_h - u\| \| y \|.
$$
It follows that
$$
\limsup_{t \to + \infty} | \langle \sigma_t x  - u, y \rangle | \leq \left(  \frac{1}{h} \int_0^h \|S_s x- x \|
ds  \right)  \|y \|  + \| u_h - u\| \| y \|
$$
for all $h>0$. Letting $h \to 0$ we obtain $\sigma_t x \wto u$.\bx

\noi Set $\bar z_n=\sfrac{1}{\sigma_n}\sum_{k=1}^n\lambda_kz_k$, where $z_n$ is given in \eqref{euler}. A general
result on convergence in average is the following from \cite[Bruck]{Bru2}:

\begin{thm}\label{Eulav}
Assume $\sum \| z_n- z_{n-1}\|^2< \infty$. If $\S \neq \emptyset$, then $z_n$ converges weakly in average to $w=\limn
P_{\S} z_n$. Otherwise $\limn\|\bar z_n\|=\infty$.
\end{thm}

\dem We first prove that  $\Omega[\bar
z_n] \subset  \S$ which is  {\bf (b')}. Then we show,  if $\S$ is non empty, that the sequence of  projections $\zeta_n=P_{\S} z_n$ converge strongly  to some $\zeta \in \S$ which is {\bf (a2)}
and finally that $\zeta$ is  the only weak cluster point of the bounded sequence $\{\bar z_n\}$.\\

\noi First, let $[u,v]\in A$ and set $ w_n=( z_n- z_{n+1})/\lambda_n\in A z_n$. We have
\begin{eqnarray}\label{euleraverage}
\| z_{n+1}-u\|^2 & = & \| z_n-\lambda_n  w_n-u\|^2\nonumber\\
& = & \|z_n-u\|^2+\|\lambda_n w_n\|^2+2\lambda_n\langle  w_n, u- z_n\rangle\\
& \le & \| z_n-u\|^2+\|\lambda_n w_n\|^2+2\lambda_n\langle v, u- z_n\rangle.\nonumber
\end{eqnarray}

\noi Summing up, neglecting the positive term of the telescopic sum on the left-hand side and dividing by $\sigma_n$ we
get
$$0\le \frac{\|z_1-u\|^2}{\sigma_n}+ \frac{1}{\sigma_n}\sum_{k=1}^n\| z_k- z_{k-1}\|^2+2\langle v,u-\bar z_n\rangle.$$

\noi Therefore $\liminf_{n\to\infty}\langle v,u-\bar z_n\rangle\ge 0$ and every weak cluster point of $\{\bar z_n\}$ lies in $\S$, by maximality.\\
Note that this is {\bf (b')}, hence  the counterpart of Lemma   \ref{almostcluster} and Lemma
\ref{almostclusterprox}.\\

\noi Next, take $u\in\S$. From equation \eqref{euleraverage} we get
\begin{equation}\label{euleraverage2}
\| z_{n+1}-u\|^2\le \| z_n-u\|^2+\|\lambda_n w_n\|^2.
\end{equation}
This implies the  convergence  of $\| z_{n+1}-u\|^2$  hence {\bf (a1)}  which ends the proof by using Lemma \ref{A1}.\\
So the use of the following  alternative is to identify the limit. This proof  of {\bf (a2)} parallels Lemma \ref{cvproj}.\\
Using the
parallelogram identity and the convexity of $\S$ we obtain
\begin{eqnarray*}
\|\zeta_{n+1}-\zeta_n\|^2 & = & 2\| z_{n+ 1}-\zeta_n\|^2+2\| z_{n+ 1}-\zeta_{n+ 1}\|^2-4\| z_{n+ 1}-\sfrac{1}{2}(\zeta_n+\zeta_{n+ 1})\|^2\\
& \le & 2\| z_{n+ 1}-\zeta_n\|^2-2\| z_{n+ 1}-\zeta_{n+ 1}\|^2.
\end{eqnarray*}
Inequality \eqref{euleraverage2} with $u=\zeta_n$  gives
$$
0\le\|\zeta_{n+1}-\zeta_n\|^2\le
2 \|\lambda_n w_n\|^2 +2\| z_n-\zeta_n\|^2-2\| z_{n+1}-\zeta_{n+1}\|^2.
 $$
\noi This implies that the sequence $\{\| z_n-\zeta_n\|^2 + \rho_n \}$ decreases, where
$\rho_n=\sum_{k\ge n}\|\lambda_k w_k\|^2$, which tends to $0$ as $n\to\infty$. Since
$$0\le\|\zeta_{n+p}-\zeta_n\|^2\le
2\rho_n +2\| z_n-\zeta_n\|^2-2\| z_{n+p}-\zeta_{n+p}\|^2,$$
 the sequence $\{\zeta_n\}$ is Cauchy and
converges as well  to  some  $\zeta $ in $\S$. The results now follows from Lemma \ref{A2}.
\bx

Observe that the same structure of proof could be applied to proximal sequences.\\

For a similar proof with two operators and forward-backward procedure see \cite[Passty]{Pas1}.\\

The  following  result due to  \cite[Pazy]{Paz5} of {\bf (b')} leads to a unified proof of  weak convergence in average for contractions  in the discrete   (Theorem \ref{thm:Bai_T_average})  or continuous case (Theorem  \ref{contaverage}).  Note that the first step assumes  $\S \neq \emptyset$
and then one  uses {\bf (a1)} to achieve the result.\\

\begin{prop}\label{unif}
Assume   $\S \neq \emptyset$, then $\Omega [\sigma _t x] \subset \S$.
\end{prop}

\dem For $t,h\ge 0$ we have
\begin{eqnarray*}
0  &  \leq  &  \|S_t x - y \|^2 - \|S_{t+h} x - S_h y \|^2\\
&  =  &  \|S_tx - S_hy \|^2 - \|S_{t+h}x - S_h y \|^2 + 2 \langle S_tx - S_hy,S_hy - y \rangle + \|S_h y - y \|^2.
\end{eqnarray*}
By taking the average we deduce that
$$
0 \leq \frac{1}{t}  \int_0^t   [\|S_s x - S_h y \|^2 - \|S_{s+h} x - S_h y \|^2 ] ds + 2 \langle  \sigma _t x - S_h y,
S_h y - y \rangle + \|S_h y - y \|^2.
$$
Since $\S\neq \emptyset$,  $\|S_t x - S_h y \|$ is bounded, hence letting $t \to + \infty$,   it follows that for any
$p \in \Omega[\sigma _t x]$, any $h\geq 0$ and any $y\in H$
$$
0 \leq 2 \langle  p- S_h y, S_h y - y \rangle + \|S_h y - y \|^2.
$$
Finally take $y=p$ so that $p = S_h p$, which means $p \in \S$.\bx



\section{Weak convergence}

Not all maximal monotone operators generate weakly convergent trajectories.\\

\begin{ejem}{\em Let $R:\R^2\to\R^2$ be the counterclockwise $\pi/2-$rotation and consider the evolution scheme defined by the
differential equation:
$$\dot u(t)  =  R(u(t)).$$
Note that
$\S = \{0\}$. The orbit starting at time $t=0$ from the point $u_0=r_0(\cos(\theta_0),\sin(\theta_0))$, $r>0$ is described by
$u(t)=r_0(\cos(t-\theta_0),\sin(t-\theta_0))$, which is bounded but does not have a limit as $t\to\infty$. However, the
average $\frac{1}{t}\int_0^tu(s)\ ds$ converges to $0$ as $t\to\infty$, by Theorem \ref{contaverage}.\\

Now let $x_n=r_n(\cos\theta_n,\sin\theta_n)$ satisfy
$$\frac{x_{n+1}-x_n}{\lambda_n}=R(x_{n+1}).$$
We have $r_{n+1}^2=\prod_{k=1}^n(1+\lambda_k^2)^{-1}r_0$ and $\theta_n=\theta_0+\sum_{k=1}^n\arctan(\lambda_k)$. The
sequence $r_n$ is decreasing. If $\lambda_n\notin\ltwo$ then $\limn x_n=0$; otherwise it stays bounded away from zero.
On the other hand, the argument $\theta_n$ is increasing. It converges if $\lambda_n\in\luno$ and diverges otherwise.
Observe also that $x_n$ converges in average to $0$ as $n\to\infty$, by Theorem \ref{proxav}.\\

Finally, let $z_n=\rho_n(\cos\phi_n,\sin\phi_n)$ satisfy
$$\frac{z_{n+1}-z_n}{\lambda_n}=R(z_n).$$
Here $\rho_{n+1}^2=\prod_{k=1}^n(1+\lambda_k^2)\rho_0$ and $\phi_n=\phi_0+\sum_{k=1}^n\arctan(\lambda_k)$. In this case
the sequence $r_n$ is increasing. It remains bounded and is convergent if, and only if, $\lambda_n\in\ltwo$. The
argument $\theta_n$ is increasing as well. It converges if $\lambda_n\in\luno$ and diverges otherwise. As before, $z_n$
converges in average to $0$ as $n\to\infty$, by Theorem \ref{Eulav} . }\theend\end{ejem}

\subsection*{Tools}
Assuming $\S$ non empty and using Lemma \ref{A1}, by virtue of Corollaries \ref{cor:basic_continuous_2},
\ref{lem:basic_prox_3} and \ref{cor:basic_euler}, in order to prove weak convergence of $u(t)$,  it suffices to verify
that its set of weak cluster points lie in $\S$ (condition {\bf (b)}). The key tool is the concept of {\em
demipositivity}, first developed in \cite[Bruck]{Bru1}.\\

A maximal monotone operator $A$ is {\em demipositive} if there exists $w\in \S$ such that for every sequence
$\{u_n\}\in D(A)$ converging weakly to $u$ and every bounded sequence $\{v_n\}$ such that $v_n\in Au_n$
\begin{equation}\label{Eq:demipositive}
\langle v_n,u_n-w\rangle \to 0\hskip25pt\hbox{implies}\hskip25pt u\in\S.
\end{equation}

\begin{prop}
Each of the following conditions is sufficient for a maximal monotone operator $A$ to be demipositive:
\begin{enumerate}
    \item $A=\partial\phi$, where $\phi$ is a proper lower-semicontinuous convex function having minimizers ($\S\neq\emptyset$).
    \item $A=I-T$, where $T$ is nonexpansive and has a fixed point ($\S\neq\emptyset$).
    \item The set $\S$ has nonempty interior.
    \item $A$ is odd and {\em firmly positive}, which means that there is $w\in\S$ such that $v\in Au$
            and $\langle v, u - w\rangle=0$ together imply $0\in Ax$.
    \item $A$ is firmly positive and sequentially weakly closed (its graph is sequentially weak/weak closed).
    \item $\S\neq\emptyset$ and $A$ is 3-monotone, which means that $\sum_{n=1}^3\langle y_n,x_n-x_{n-1}\rangle\ge 0$ for every set
$\{[x_n,y_n]\ |\ 1\le n\le 3\}\subset A$ ($x_0\equiv x_N$).
\end{enumerate}
\end{prop}

For demipositivity in Banach spaces see \cite[Bruck and Reich]{BrR}.\\

{\it Comments}\\
We just mention another assumption that guarantees that the weak cluster points will lie in $\S$: Let $S$ the
semi-group generated by $A$. $A$ satisfies condition $(L)$ if
$$\limty{t}\|A^0S_tx\|\le\limty{t}\left(\frac{1}{h}\|S_{t+h}x-S_tx\|\right)$$
for every $h>0$ and $x\in D(A)$. An equivalent formulation is the following: Denote by $a^0$ the element of minimal
norm in $\overline{R(A)}$. Then $A$ satisfies condition $(L)$ if, and only if, for every $x\in D(A)$ one has
$$\limty{t}A^0S_tx=a^0.$$

Unlike demipositivity, this does not impose {\em a priori} that $\S\neq\emptyset$. For instance, if $A=\partial f$ with
$f\in \Gamma_0(H)$ or if $A=I-T$ with $T$ nonexpansive, then $A$ satisfies condition $(L)$ but is not demipositive
unless $\S\neq\emptyset$.

Condition $(L)$ is essentially used in \cite [Pazy]{Paz1} to prove that the weak cluster points of the trajectory $S_tx$ lie
in $\S$. If $\S=\emptyset$ one immediately deduces that $\limty{t}\|S_tx\|=\infty$. The interested reader may find this
definition and related results in \cite[Pazy]{Paz1}.\theend

\subsection{Continuous dynamics}

\noi The following classical result of weak convergence for demipositive operators was proved in \cite[Bruck]{Bru1}.

\begin{thm}\label{demiweak}
If $A$ is demipositive then $u(t)$ converges weakly as $t\to\infty$ to an element of $\S$.
\end{thm}

\dem By Corollary \ref{cor:basic_continuous_2} and Opial's Lemma it suffices to prove $\Omega[u(t)] \subset \S$, which
is {\bf (b)}. Let $w\in\S$ satisfy \eqref{Eq:demipositive} and let $u(t_n)  \wto u$ as $n\to\infty$.  The sequence
$\dot u(t_n)$ is bounded by Theorem \ref{exuni}.  Let $\theta_w(t)=\frac{1}{2}\|u(t)-w\|^2$, thus
$\dot\theta_{w}(t)=\langle \dot u(t_n),u(t_n)-w\rangle$. Since $\theta_{w}$ is bounded by Corollary
\ref{cor:basic_continuous_2}, $\dot\theta_{w}\in L^1$ and there is a subsequence $t_{n_k}$ of $t_n$ such that
$\dot\theta_{w}(t_{n_k})\to 0$ as $k\to\infty$. But $u(t_{n_k})\wto u$, so $u\in\S$ by demipositivity.\bx

{\it Comments}\\
Theorem \ref{demiweak} was extended in \cite[Passty]{Pas} to the class of $\varphi$-{\em demipositive}
operators.\theend

\subsection{Proximal sequences}

A first detailed study of the asymptotic behavior of the  proximal sequence $\{x_n\}$ was performed in
\cite[Rockafellar]{Roc}, when the stepsizes are bounded away from zero. The author also considers an inexact version of
the algorithm. The next convergence results under more general hypotheses are investigated in \cite[Br\'ezis and Lions]{BrL}.\\

Recall that $\sigma_n=\sum_{m\le n}\lambda_m$ and $\tau_n=\sum_{m\le n}\lambda_m^2$.

\begin{thm}\label{weakprox} Assume $\S\neq\emptyset$. If $\{\lambda_n\}\notin\ltwo$
then $x_n$ converges weakly to some $x^*\in \S$. Moreover, $\|y_n\|\leq d(x_0,\S)\tau_n^{-1/2}$.
\end{thm}
\dem By Lemmas \ref{lem:basic_prox_1} and \ref{lem:basic_prox_2}, we have for any $x \in \S$
$$
\| y_n\|^2\tau_n\le \sum_{k\le n} \lambda_k^2 \| y_k\|^ 2 \leq \| x_0 - x \|^2.
$$
 $\tau_n\to \infty$ implies $\|y_n\|\to 0$.
  Since $-y_n \in Ax_n$, we deduce that $\Omega[x_n] \subset \S$, which is {\bf (b)}, by
 Proposition \ref{prop:ws_sw}.
 We conclude by Corollary \ref{lem:basic_prox_3} and Opial's Lemma \ref{Op}.\bx

The following result, adding the demipositivity hypothesis, is also from \cite[Br\'ezis and Lions]{BrL}:

\begin{thm}\label{weakprox_demi}
If $A$ is demipositive then $x_n$ converges weakly to some $x^*\in \S$.
\end{thm}

\dem As above, using Corollary \ref{lem:basic_prox_3}  the result follows from Opial's Lemma \ref{Op}  if $\Omega[x_n]
\subset \S$ which is {\bf (b)}. Let $x_{n_k} \wto x$ and $w$ be the element in $\S$  used  in the definition of demipositivity
\eqref{Eq:demipositive}. Using Lemma \ref{lemweakprox} below we construct another subsequence $\{x_{m_k}\}$ such that
both  $\|x_{m_k} - x_{n_k}\|$ and $ \langle x_{m_k} -w, y_{m_k}\rangle$ tend to $0$ as $k\to\infty$. Since $x_{m_k}\wto
x$ and $A$ is demipositive, $x$ must belong to $\S$. \bx

\begin{lem}\label{lemweakprox}
 Let $\{x_n\}$ be a proximal sequence and $w\in\S$. For each
$\varepsilon>0$,  there is $N$ such that:  for any $n \geq N$, there exists $m\in\N$ satisfying $N \leq m\leq n$, $\|x_m -x_n \| \leq\varepsilon$ and $\langle -y_m, x_m - w \rangle \leq \varepsilon$.
\end{lem}

\dem For each $w\in \S$ we have $\| x_{k-1} - w \| ^2 \geq \| x_k - w \|^2 + 2 \lambda _k \langle -y_k , x_k - w
\rangle$ and so
\begin{equation}\label{brl_weak}
\sum_k \lambda _k \langle y_k , w -x_k \rangle < \infty
\end{equation}
where all terms are nonnegative by monotonicity. Given $\eps>0$, define $P= \{ k\in\N\ |\ \langle y_k , w-x_k \rangle
\geq \varepsilon \}$ so that $\sum _{ k\in P} \lambda _k < \infty$. Since $\|x_{k-1} - x_k\| = \lambda_k \|y_k\|$,
Lemma \ref{lem:basic_prox_1} implies $\sum _{k\in
P} \|x_{k-1} - x_k\|  < \infty$.\\
 Let $N_1$ so that $ \sum _{k\in P, k\ge N_1} \|x_{k-1} - x_k\| < \varepsilon$. By virtue of \eqref{brl_weak},
since $\{\lambda_n\}\notin\luno$ there is $N\geq N_1$ with $\langle y_N , w-x_N \rangle \leq \varepsilon$. Consider
$n\geq N$: if $n\notin P$ we choose $m=n$. If $n\in P$, let $m=\max\{k<n\ |\ k\notin P\}$. Since $m\geq N_1$ and all
integers
between $m$ and $n$ are in $P$, we have $\|x_m - x_n \| \leq \sum _{m<k\leq n}\|x_{k-1} - x_k\| \leq  \varepsilon$. \bx\\

{\it Comments}
\begin{enumerate}
        \item Theorem \ref{weakprox_demi} is still true if the sequence satisfies
            $ \|x_n - (I + \lambda_n A) ^{-1} x _{n-1} \| \leq \varepsilon_n $ with $\sum \varepsilon_n < \infty$.
            This is proved in \cite[Br\'ezis and Lions]{BrL} and can also be derived using asymptotic equivalence results
            in Section \ref{section:almost-equivalence} (see \cite{AlP1,AlP2}).
        \item In uniformly convex Banach spaces with Fr\'echet differentiable norm there
            is weak convergence in the following cases (see \cite[Reich]{Rei3}):
        \begin{enumerate}
                \item $\{\lambda_n\}$ does not converge to zero, or
                \item The modulus of convexity of the space satisfies $\delta(\eps)\ge K\eps^p$ for some $K>0$ and
                    $p\ge 2$ and $\sum\lambda_n^p=\infty$.
        \end{enumerate}
    \item Demipositive can be replaced by $\varphi$-demipositive (see \cite[Passty]{Pas}).\theend
\end{enumerate}

\subsection{Euler sequences}
Let $\{z_n\}$ be an Euler sequence and recall that $w_n = \frac{z_{n+1} - z_n}{\lambda_n}  \in - Az_n$.

\begin{thm}\label{thm:weak_euler}
Let $A$ be demipositive and assume $\{\lambda_n\}\in\ltwo$ and $\{w_n\}$ bounded. Then $z_n$ converges weakly to some $z \in \S$.
\end{thm}

\dem If $y\in\S$,  Corollary \ref{cor:basic_euler}  shows that the sequence $\|z_n-y\|$ is convergent. On the other hand, equality
\eqref{Eq:Euler_1} implies $\sum_{n\ge 1}\lambda_n\langle w_n,y-z_n\rangle<\infty$. One concludes as in Theorem
\ref{weakprox_demi} using an analogue of Lemma \ref{lemweakprox}.\bx

{\it Comments}\\
The previous result from \cite[Bruck and Reich]{BrR} works for demipositive operators in ``a few" Banach spaces, namely
$X=L^{2m}$, $m\in\N$ or $X=\ell^p$, $p\in(1,\infty)$. \theend

A related result from \cite[Reich]{Rei3} is the following (and  holds in uniformly convex Banach spaces with
Fr\'echet-differentiable norm):

\begin{prop} Let $T$ be non-expansive, $A=I-T$ and $\{\lambda_n\}$ satisfying $0\le\lambda_n\le 1$ and
$\sum\lambda_n(1-\lambda_n)=\infty$. If $\S \not= \emptyset$
then $\{z_n\}$ converges weakly to a  point in $\S$.
\end{prop}

If $A=\partial f$ with $f\in \Gamma_0(H)$ and $\hbox{dim}(H)<\infty$ one can circumvent the difficulties of Lemma
\ref{lemweakprox} and provide a simpler proof of Theorem \ref{thm:weak_euler}. Let
$$z_{n+1} \in z_n - \lambda_n \partial f (z_n)$$
and $w_n$ as above.

\begin{thm}
Assume $\S\neq\emptyset$ and $\hbox{dim}(H)<\infty$. If $\sum\|z_n-z_{n-1}\|^2<\infty$ then $z_n$ converges to a
minimizer of $f$.
\end{thm}

\dem Lemma  \ref{lem:Euler_liminf} gives $\liminf\limits_{n\to\infty}f(z_n)=f^*$. Since $\{z_n\}$ is bounded and
the space is finite dimensional, there is a subsequence $\{z_{n_k}\}$ such that
$\lim\limits_{k\to\infty}f(z_{n_k})=f^*$ and $\lim\limits_{k\to\infty}\|z_{n_k}-z\|=0$ for some $z\in H$. Since $z$
must be in $\S$ by lower-semicontinuity, Corollary \ref{cor:basic_euler} implies
$\lim\limits_{n\to\infty}\|z_{n}-z\|=0$, which means $z_n$ converges to $z$.\bx

The preceding result from \cite{She} was pointed out to the authors by R. Cominetti.


\section{Strong convergence}

Even if $A=\partial\phi$ with $\phi\in\Gamma_0(H)$ having minimizers, the trajectory $u(t)$ need not converge strongly
as $t\to\infty$. This is shown by Baillon's example in \cite[Baillon]{Bai}: the author defines a function
$\phi\in\Gamma_0(\ltwo)$ having minimizers and proves that the trajectories converge weakly but not strongly.\\

This also true for the proximal point algorithm. Even if $A=\partial\phi$ with $\phi\in\Gamma_0(H)$ having minimizers,
a sequence satisfying \eqref{prox} need not converge strongly. This was proved in \cite[G\"uler]{Gul} using Baillon's
example and the equivalence techniques from \cite[Passty]{Pas}. A simpler example of this type can be found in
\cite{BBDHV} and can be retranslated to provide a new counterexample for strong convergence of the continuous
trajectory, different from that of Baillon.

\subsection*{Conditions}

We introduce here a series of conditions, mainly of geometric nature, that will be used to obtain strong convergence of the process in the  continuous or discrete set-up.  \\

\noi{\bf Strong monotonicity.} Let $\alpha>0$. An operator $A$ is {\em $\alpha$-strongly monotone} if for all
$[x,x^*],[y,y^*]\in A$ one has
$$\langle x^*-y^*,x-y\rangle\ge\alpha\|x-y\|^2.$$
Observe that if $A$ is strongly monotone and $Ax\cap Ay\neq\emptyset$, then $x=y$. If $A$ is $\alpha$-strongly monotone
then  $J_{1/\alpha}^A$ is a strict contraction. Therefore it has a fixed point $p$ and only one, say $\S=\{p\}$.
Strongly monotone operators are demipositive.\\

\noi Clearly, if $A$ is monotone, then $A+\alpha I$ is $\alpha$-strongly monotone. Also, subdifferentials of proper,
lower-semicontinuous strongly convex functions are strongly monotone.\\

A weaker notion of strong monotonicity found for instance in \cite[Pazy]{Paz1} is the following: $A$ is
$\alpha$-strongly monotone if $\S\neq\emptyset$ and
$$\langle A^0x,x-P_{\S}x\rangle\ge \alpha\|x-P_{\S}x\|^2$$
for every $x\in D(A)$. In this case the set $\S$ need not be a singleton. Proposition \ref{prop:strong_continuous}
below also holds if $A$ is strongly monotone in this sense but the proof is more involved.\\

\noi{\bf Solution set $\S$ with nonempty interior.} If $p\in\interior\S$ then there is $r>0$ such that the ball
$B(p,r)$ of radius $r$ centered at $p$ is contained in $\S$. Then $\langle u^*,u-p+rh\rangle\ge 0$ for all $[u,u^*]\in
A$ and all $h\in H$ with $\|h\|\le 1$. Therefore $\langle u^*,u-p\rangle\ge r\langle u^*,-h\rangle$ and
\begin{equation}\label{Eq:nonempty}
r\|u^*\|=r\sup_{\|h\|\le 1}\langle u^*,-h\rangle\le\langle u^*,u-p\rangle.
\end{equation}

\noi{\bf The NR convergence condition.}\label{Nev-Rei} A maximal monotone operator $A$ on $H$ satisfies the {\em NR
convergence condition} if $\S\neq\emptyset$ and for every bounded sequence $[x_n,y_n]\in A$ one has
$$\liminf\limits_{n\to\infty}\langle y_n,x_n-P_{\S}x_n\rangle=0
\qquad\hbox{implies}\qquad\liminf\limits_{n\to\infty}\|x_n-P_{\S}x_n\|=0.$$

Strongly monotone operators satisfy this condition. So do operators having compact resolvent (see below) and those
satisfying $\langle y,x-P_{\S}x\rangle>0$ for all $[x,y]\in A$ such that $x\notin\S$.\\

The NR convergence condition can be easily stated in a Banach space $X$ by means of
the duality mapping. The results below hold when both $X$ and $X^*$ are uniformly convex.
The interested reader can consult \cite[Nevanlinna and Reich]{NeR} and \cite[Bruck and Reich]{BrR}.\\

\noi{\bf Compactness.} The strong $\omega$-limit set of a trajectory $u:[0,\infty)\to H$ is the set
$\omega[u(t)]=\bigcap\limits_{t>0}\overline{\{u(s) :s\ge t\}}$. For a sequence $\{x_n\}$ it is defined by
$\omega[x_n]=\bigcap\limits_{n\in\N}\overline{\{x_k :k\ge n\}}$.\\

By virtue Lemma \ref{A1} the sets $\omega[u(t)]\cap\S$ and $\omega[x_n]\cap\S$ contain, at most, one element.\\

If $\S\neq\emptyset$ and $J_1^A$ is a compact operator (maps bounded sets to relatively compact sets) then
$\omega[u(t)]\neq\emptyset$ for every trajectory $u$ satisfying \eqref{continuous} (see Theorem 11.8 in
\cite[Pazy]{Paz1}) and $\omega[x_n]\neq\emptyset$ for every sequence $\{x_n\}$ satisfying \eqref{prox}. \\

For instance, if $A=\partial f$ and the set $\{\ u\in H\ |\f(u)+\|u\|^2\le M\ \}$ is compact for each $M\ge 0$, then
$J_1^A$ is compact. This case was first studied in \cite[Br\'ezis]{Bre}.\\

\noi{\bf Symmetry.} An operator $A$ is odd if $D(A)=-D(A)$ and $A(-x)=-Ax$ for all $x\in D(A)$.\footnote{A weaker
notion is that $A^0(-x)=-A^0x$. The results below still hold but the proofs become more technical.} This is the case,
for instance, if $A=\partial f$ and $f$ is even. If $A$ is odd, the semigroup generated is odd as well (see, for
instance, \cite[Pazy]{Paz1}). On the other hand, it is easy to see that $J_\lambda^A$ is odd for each $\lambda>0$ if
$A$ is odd.\\

Notice also that if $A$ is odd then $\S\neq\emptyset$. Moreover, $0\in\S$. To see this, take $x\in D(A)$ and let
$[x,y],[-x,-y]\in A$. We have
\begin{eqnarray*}
4\langle y-0,x-0\rangle & = & \langle y+y,x+x\rangle\\
& = & \langle y-(-y),x-(-x)\rangle\\
& \ge & 0.
\end{eqnarray*}
Then $0\in A0$ by Lemma \ref{lem:mm_sufficient}.\\

\noindent{\bf Asymptotic regularity.} A trajectory $u$ is {\em asymptotically regular} if $\limty{t}\|u(t+h)-u(t)\|=0$
for each $h\ge 0$. A sequence $\{x_n\}$ is asymptotically regular if $\limty{n}\|x_{n+m}-x_n\|=0$ for each $m\in\N$.\\

\noi{\it Comments}\\
Recall that the notion of weak asymptotic regularity was mentioned in Proposition \ref{PP3} as a characterization of
weak convergence of the trajectories satisfying \eqref{continuous}.\theend

\subsection{Continuous dynamics}

\noindent{\bf Strong monotonicity.}

\begin{prop}\label{prop:strong_continuous}
If $A$ is $\alpha$-strongly monotone for some $\alpha>0$ then $u(t)$ converges strongly to the unique $p\in\S$ as
$t\to\infty$.
\end{prop}

\dem Strong monotonicity implies
$$\frac{1}{2}\frac{d}{dt}\|u(t)-p\|^2 = \langle \dot u(t),u(t)-v(t)\rangle \le -\alpha \|u(t)-p\|^2$$
and so $\|u(t)-p\|\le e^{-2\alpha t}\|u_0-p\|$.\bx

{\it Comments}\\
The previous result can be extended in the following way: Let $X$ be a Banach space such that $X$ and $X^*$ are
uniformly convex. In \cite[Nevanlinna and Reich]{NeR} the authors prove that if $A$ satisfies NR convergence condition
then $u(t)$ converges strongly to a point in $\S$ as $t\to\infty$. If only $X^*$ is uniformly convex, the result
remains true provided $Ax$ is proximinal and convex for every $x$ (see \cite[Bruck and Reich]{BrR}). If neither $X$ nor
$X^*$ is uniformly convex, the result is still true if the semigroup is differentiable (see \cite[Nevanlinna and
Reich]{NeR}). \theend

\noindent{\bf Solution set with nonempty interior.}

\begin{prop}\label{stronginterior}
Assume $\interior\ \S\neq\emptyset$. Then $u(t)$ converges strongly as $t\to\infty$ to a point in $\S$.
\end{prop}

\dem If $B(p,r)\subset\S$, inequality \eqref{Eq:nonempty} implies
\begin{eqnarray*}
r\|u(t)-u(s)\| & \le & r\int_s^t\|\dot u(\tau)\|\ d\tau\\
& \le & -\int_s^t\langle \dot u(\tau),u(\tau)-p\rangle\ d\tau\\
& = & \frac{1}{2}\|u(s)-p\|^2-\frac{1}{2}\|u(t)-p\|^2.
\end{eqnarray*}
Since $\|u(t)-p\|$ is convergent by Corollary \ref{cor:basic_continuous_2}, $u(t)$ has the Cauchy property.\bx

{\it Comments}\\
Theorem 4 in \cite[Nevanlinna and Reich]{NeR} shows that this result remains true if $X$ and $X^*$ are uniformly
convex. In the same paper, the authors give a counterexample in $\C([0,1];\R)$. See also \cite[Bruck and Reich]{BrR}.
\theend

\noindent{\bf Compactness.}

\begin{prop}\label{omegalimit}
If $\omega[u(t)]\cap\S\neq\emptyset$ then $u(t)$ converges strongly to some $p\in\S$.
\end{prop}

\dem If $p\in\omega[u(t)]\cap\S$ then $\|u(t)-p\|$ is decreasing and $\liminf\limits_{t\to\infty}\|u(t)-p\|=0$. Hence
$u(t)\to p$ as $t\to\infty$.\bx

{\it Comments}\\
If $\S$ has nonempty interior then $A$ is demipositive and $\omega[u(t)]\neq\emptyset$ for every trajectory $u$
satisfying \eqref{continuous}. Every strong cluster point is also a weak cluster point, that must lie in $\S$ by
demipositivity. Hence $\omega[u(t)]\cap\S\neq\emptyset$ and Proposition \ref{stronginterior} can also be deduced from
Proposition \ref{omegalimit}. \theend

\noi{\bf Symmetry.}

\begin{prop}
 If $A=\partial f$ and $f\in \Gamma_0(H)$ is even then $u(t)$ converges strongly as
$t\to\infty$ to a point in $\S$.
\end{prop}

\dem Take $s>0$ and define $\gamma(t)=\|u(t)\|^2-\|u(s)\|^2-\hbox{$\frac{1}{2}$}\|u(t)-u(s)\|^2$. For $t\in[0,s]$ one
has
$$\dot \gamma(t)=\langle
\dot u(t),u(t)+u(s)\rangle\le f(-u(s))-f(u(t))=f(u(s))-f(u(t))\le 0.$$ Therefore, $\gamma(t)\ge\gamma(s)=0$ and so
$$\frac{1}{2}\|u(t)-u(s)\|^2\le\|u(t)\|^2-\|u(s)\|^2.$$
Since $0\in\argmin\ f$, $\|u(t)\|$ converges as $t\to\infty$ so $u(t)$ has the Cauchy property.\bx

For general $A$ one has to assume additional hypotheses on the trajectory:

\begin{prop}\label{prop:Pazy_BBR}
Let $A$ be odd. If $u$ is asymptotically regular
 then $u(t)$ converges strongly to
some $p\in\S$ as $t\to\infty$.
\end{prop}

\dem Let us use the semigroup notation $u(t)=S_tx$. If $A$ is odd then $0\in\S$ and
\begin{eqnarray*}
\|S_{t+h+s}x+S_{t+s}x\| & = & \|S_{t+h+s}x-S_{t+s}(-x)\|\\
& \le & \|S_{t+h}x-S_t(-x)\|\\
& = & \|S_{t+h}x+S_tx\|
\end{eqnarray*}
for each $h\ge 0$ so that
\begin{equation}\label{Eq:odd_continuous}
\limty{t}\|S_tx+S_{t+h}x\|\le\|S_tx+S_{t+h}x\|.
\end{equation}
Since $0\in\S$ the limit $d=\limty{t}\|S_tx\|$ exists. Moreover, the fact that
$\|2S_tx\|\le\|S_tx+S_{t+h}x\|+\|S_tx-S_{t+h}x\|$ implies
$$2d\le\limty{t}\|S_tx+S_{t+h}x\|\le\|S_tx+S_{t+h}x\|$$ for each $t,h$ by asymptotic regularity and inequality \eqref{Eq:odd_continuous}. Finally,
\begin{eqnarray*}
\|S_{t+h}x-S_tx\|^2 & = & 2\|S_tx\|^2+2\|S_{t+h}x\|^2-\|S_{t+s}x+S_tx\|^2\\
& \le & 4\|S_tx\|^2-4d^2
\end{eqnarray*}
and so $\{S_tx\}$ has the Cauchy property. \bx

\noi{\it Comments}\\
Without the asymptotic regularity assumption, strong convergence holds for the averages when $S$ is odd, as proved in
(see \cite[Baillon]{Bai3}).\theend

\subsection{Proximal sequences}

\noindent{\bf Strong monotonicity.}

\begin{prop}
If $A$ is $\alpha$-strongly monotone for some $\alpha>0$ then $x_n$ converges strongly to the unique $p\in\S$ as $n \to \infty$.
\end{prop}

\dem Strong monotonicity implies
\begin{eqnarray*}
\alpha\lambda_n\|x_n-p\|^2 & \le & \langle x_{n-1}-x_n,x_n-p\rangle\\
& = & \langle x_{n-1}-p,x_n-p\rangle-\|x_n-p\|^2\\
& \le & \|x_n-p\|\left(\|x_{n-1}-p\|-\|x_n-p\|\right)
\end{eqnarray*}
so that
$$\alpha\sum_{n=1}^\infty\lambda_n\|x_n-p\|\le\|x_0-p\|<\infty.$$
Since the sequence  $\|x_n-p\|$ is decreasing this implies $\limn\|x_n-p\|=0$.\bx

\noindent{\bf Solution set with nonempty interior.}

\begin{prop}
Let $A$ be maximal monotone with $\interior\ \S\neq\emptyset$. Then $x_n$ converges strongly as $n\to\infty$.
\end{prop}

\dem If $B(p,r)\subset\S$ inequality \ref{Eq:nonempty} gives $r\|x_{k-1}-x_k\|\le \langle x_{k-1}-x_k,x_k-p\rangle$ and
so
\begin{eqnarray*}
r\|x_{k-1}-x_k\| & \le & \langle x_{k-1}-p,x_k-p\rangle-\|x_k-p\|^2\\
& \le & \|x_0-p\|\left(\|x_{k-1}-p\|\frac{}{}\!-\|x_k-p\|\right)
\end{eqnarray*}
by Corollary \ref{lem:basic_prox_3}. Hence
\begin{eqnarray*}
r\ \|x_n-x_m\| & \le & r\sum_{k=n+1}^m\|x_{k-1}-x_k\|\\
& \le & \|x_0-p\|\left(\|x_n-p\|\frac{}{}-\|x_m-p\|\right).
\end{eqnarray*}
Since $\|x_n-p\|$ is convergent, $x_n$ is a Cauchy sequence.\bx

\noi{\bf The NR convergence condition.}

\noi A fairly general result is the following, from \cite[Nevanlinna and Reich]{NeR}:

\begin{thm}
If $A$ satisfies the NR convergence condition then $x_n$ converges strongly as $n\to\infty$.
\end{thm}

\dem Setting $j_n=x_n-P_{\S}x_n$ we have
\begin{eqnarray*}
\|j_n\|^2+\lambda_{n}\langle y_{n},j_{n}\rangle & = & \langle x_{n-1}-P_{\S}x_{n},j_{n}\rangle\\
& = & \langle j_{n-1},j_{n}\rangle + \langle P_{\S}x_{n-1}-P_{\S}x_{n},x_n-P_{\S}x_n\rangle\\
& \le & \|j_{n-1}\|\ \|j_n\|\\
& \le & \frac{1}{2}\left[\|j_{n-1}\|^2+\|j_n\|^2\right].
\end{eqnarray*}
Thus $\|j_n\|^2+2\lambda_{n}\langle y_{n},j_{n}\rangle\le \|j_{n-1}\|^2$ and $\sum_{n=1}^\infty\lambda_{n}\langle
y_{n},j_{n}\rangle<\infty$. Since $\langle y_{n},j_{n}\rangle\ge 0$ one must have $\liminf_{n\to\infty}\langle
y_n,j_n\rangle=0$. The sequences $\{x_n\}$ and $\{y_n\}$ are bounded, and the convergence condition implies
$\liminf_{n\to\infty}\|x_n-P_{\S}x_n\|=0$. Since $\|x_n-P_{\S}x_n\|$ is nonincreasing, it must converge to $0$. On the
other hand, the sequence $\|x_n-p\|$ is nonincreasing for each $p\in\S$. In particular,
$\|x_{n+m}-P_{\S}x_n\|\le\|x_n-P_{\S}x_n\|$ and therefore $\|x_{n+m}-x_n\|\le 2\|x_n-P_{\S}x_n\|$. We conclude that
$x_n$ converges strongly to some $p\in\S$ as $n\to\infty$.\bx

\noi{\bf Compactness.}

\begin{prop}
If $\omega[x_n]\cap\S\neq\emptyset$ then $x_n$ converges strongly to some $p\in\S$.
\end{prop}

\dem If $p\in\omega[x_n]\cap\S$ then $\|x_n-p\|$ is decreasing and $\liminf\limits_{n\to\infty}\|x_n-p\|=0$. \bx

\noi{\bf Symmetry.}

For even functions we have the following result from \cite[Br\'ezis and Lions]{BrL}:

\begin{prop}
If $A$ is the subdifferential of an even function in $f\in\Gamma_0(H)$ then $x_n$ converges strongly as $n\to\infty$.
\end{prop}

\dem Recall that $2\lambda_n(f(u)-f(x_n))\ge\|u-x_n\|^2-\|u-x_{n-1}\|^2$. Let $m\ge n$ and take $u=-x_m$. Since
$n\mapsto f(x_n)$ is decreasing we have $\|x_m+x_n\|\le\|x_m+x_{n-1}\|$ and the function $n\mapsto\|x_m+x_n\|$ is
decreasing. In particular $\|x_m+x_m\|\le\|x_m+x_n\|$, thus $4\|x_m\|^2\le\|x_m+x_n\|^2$. We have
$2\|x_n\|^2+2\|x_m\|^2=\|x_m+x_n\|^2+\|x_m-x_n\|^2\ge 4\|x_m\|^2+\|x_m-x_n\|^2$, so that $\|x_m-x_n\|^2\le
2\|x_n\|^2-2\|x_m\|^2$. Since $\|x_n\|$ converges as $n\to\infty$ this proves that $x_n$ is a Cauchy sequence.\bx

As before, asymptotic regularity is required for a general $A$:

\begin{prop}\label{prox_odd}
Let $A$ be odd. If $\{x_n\}$ is asymptotically regular then $x_n$ converges strongly to some
$p\in\S$ as $n\to\infty$.
\end{prop}

\dem First, one easily verifies that $0\in\S$ and that the sequence $\|x_{n+k}+x_n\|$ is decreasing for each $k\in\N$.
Finally one concludes as in the proof of Proposition \ref{prop:Pazy_BBR}. \bx

{\it Comments}\\
\noi Without asymptotic regularity on can still prove strong convergence of the averages (see \cite[Lions]{Lio}) for
odd operators. This was first proved in \cite[Baillon]{Bai33} in the case $\lambda_n\equiv\lambda$. \theend

\subsection{Euler sequences}

\noi{\bf Strong monotonicity.}

\begin{prop}
Let $A$ be $\alpha$-strongly monotone. If $\sum\|z_n-z_{n-1}\|^2<\infty$ then $z_n$ converges strongly to the unique
$p\in\S$ as $n \to \infty$.
\end{prop}

\dem Strong monotonicity implies
$$2\alpha\lambda_n\|z_n-p\|^2 + \|z_{n+1}-p\|^2 \le \|z_n-p\|^2 + \lambda_n^2\|w_n\|^2.$$
Therefore
$$2\alpha\sum_{n=1}^\infty\lambda_n\|z_n-p\|^2\le\|z_0-p\|^2+\sum\lambda_n^2\|w_n\|^2<\infty.$$
This implies $\liminf\limits_{n\to\infty}\|x_n-p\|=0$. But $\|x_n-p\|$ converges by Corollary \ref{cor:basic_euler}.\bx

\noi{\bf Solution set with nonempty interior.}

\begin{prop}
Assume $\interior\ \S\neq\emptyset$. If $\sum\|z_n-z_{n-1}\|^2<\infty$ then $z_n$ converges strongly as $n\to\infty$.
\end{prop}

\dem If $B(p,r)\subset\S$ inequalities \eqref{Eq:nonempty} and \eqref{monEul} together give
$$2r\lambda_n\|w_n\| + \|z_{n+1}-p\|^2 \le \|z_n-p\|^2 + \lambda_n^2\|w_n\|^2.$$
This implies the sequence $\lambda_n\|w_n\|=\|z_{n+1}-z_n\|$ is in $\luno$ and so $z_n$ converges.\bx

\noi{\bf The NR convergence condition.}

\begin{thm}\label{nereuler}
Assume $\sum\|z_{n+1}-z_n\|^2 <\infty$ and $w_n$ is bounded. If $A$ satisfies the NR convergence condition then
$\{z_n\}$ converges strongly as $n\to\infty$.
\end{thm}

\dem To simplify notation write $j_n=z_n-P_S z_n$. We have
$$\|j_{n+1}\|^2 \le \|z_{n+1}-P_\S z_n\|^2 = \|j_n+\lambda_nw_n\|^2
= \|j_n\|^2-2\lambda_n\langle -w_n,j_n\rangle+\lambda_n^2\|w_n\|^2.$$

By hypothesis and Corollary \ref{cor:basic_euler} the sequence $[z_n,-w_n]$ is bounded. Moreover,
$$\sum\limits_{n=1}^\infty\lambda_n\langle -w_n,j_n\rangle <\infty.$$
But $\langle -w_n,j_n\rangle\ge 0$ and so $\liminf\limits_{n\to\infty}\langle w_n,j_n\rangle=0$ and the convergence
condition implies $\liminf\limits_{n\to\infty}\|j_n\|=0$. This sequence being convergent we have $\limn j_n=0$.
Finally, $\|z_{n+m}-z_n\|\le 2\|j_n\|$ and so $z_n$ converges as $n\to\infty$.\bx

{\it Comments}\\
The previous result holds if $X$ and $X^*$ are uniformly convex (see \cite[Nevanlinna and Reich]{NeR}).\\

According to \cite[Bruck and Reich]{BrR}, the convergence condition can be replaced by $\interior\ \S\neq\emptyset$. In
that case, if $X$ is not uniformly convex it suffices that $Ax$ be proximinal and convex for each $x$. On the other
hand, according to \cite[Nevanlinna and Reich]{NeR}, the conclusion of Theorem \ref{nereuler} is still true, even if
$X$ and $X^*$ are not uniformly convex, provided $\S$ is proximinal and $A$ is {\em accretive in the sense of Browder}.
\theend

\noi{\bf Compactness.}

\begin{prop}
Assume that $\sum\|z_{n+1}-z_n\|^2<\infty$ and $\omega[z_n]\cap\S\neq\emptyset$. Then $z_n$ converges strongly to some
$p\in\S$.
\end{prop}

\dem The argument is the same as in Proposition \ref{omegalimit} by virtue of Corollary \ref{cor:basic_euler}.\bx

\noi{\bf Symmetry.}\\

The following results uses the same ideas as in Propositions \ref{prop:Pazy_BBR} and \ref{prox_odd} but is apparently
new:

\begin{prop}
Let $T$ be non-expansive, $A=I-T$ and $\lambda_n\equiv 1$ so that $z_n=T^nz_0$. If $T$ is odd and
$\sum\|z_{n+1}-z_n\|^2 <\infty$ then $\{z_n\}$ is strongly convergent.
\end{prop}

\dem Since $T$ is odd one easily deduces that the sequence $\|z_{n+k}+z_n\|$ is decreasing for each $k$. From the fact
that $\sum\|z_{n+1}-z_n\|^2 <\infty$ we can draw two conclusions: In the first place, Corollary \ref{cor:basic_euler}
implies $d=\limn\|z_n\|$ exists because $0\in\S$. On the other hand, the sequence $z_n$ is asymptotically regular, so
$\limn\|z_n-z_{n+k}\|$ exists for each $k$. As a consequence, $2d\le\|z_{n+k}+z_n\|$ for each $n$ and $k$. One
concludes as in the proof of Proposition \ref{prop:Pazy_BBR}.\bx

{\it Comments}\\
Without any further assumptions, $T^nz$ converges strongly in average if $T$ (see \cite[Baillon]{Bai33}).\theend


\section{Asymptotic equivalence}\label{section:almost-equivalence}

In this section we explain how to deduce qualitative information on the asymptotic behavior of the systems defined by
 \eqref{continuous}, \eqref{prox} and \eqref{euler}. We provide a comparison tool that guarantees that two evolution
systems share certain asymptotic properties. For the complete abstract theory see \cite[Alvarez and Peypouquet]{AlP2}.\\

\subsection{Evolution systems}

\noi Let $C$ be a convex subset of a Banach space $X$ and let $I$ denote the identity operator in $X$. An {\em
evolution system (ES)} on $C$ is a family $\{V(t,s)\ :\ t\ge s\ge 0\ \}$ of maps from $C$ into itself satisfying:
    \begin{itemize}
        \item [i)] $V(t,t)=I$; and
        \item [ii)] $V(t,s)V(s,r)=V(t,r)$.
    \end{itemize}
Let $L>0$. An evolution system is $L$-{\em Lipschitz} if it satisfies
    \begin{itemize}
        \item [iii)] $\|V(t,s)x-V(t,s)y\|\le L\|x-y\|$
    \end{itemize}
and is {\em contracting (CES)} if it is $1$-Lipschitz.


\begin{ejem}\label{E:es1}{\em Let $F$ be a (possibly multivalued) function from $[t_0,\infty)\times
C$ to $C$. Suppose that for every $s\ge t_0$ and $x\in C$ the differential inclusion
$u'(t) \in F(t,u(t))$, with initial condition $u(s) = x$, has a unique solution $u_{s,x}:[s,\infty)\mapsto C$. The
family $U$ defined by $U(t,s)x=u_{s,x}(t)$ is an evolution system on $C$. If $X$ is Hilbert space and $F(t,x)=-A_tx$,
where $\{A_t\}$ is a family of maximal monotone operators, then the corresponding $U$ is a \ces. }\wbx
\end{ejem}

\begin{ejem}{\em \label{E:es2} Take a strictly increasing unbounded sequence $\{\sigma_n\}$ of
positive numbers and set $\nu(t)=\max\{n\in\N\ |\ \sigma_n\le t\}$. Consider a family $\{F_n\}$ of functions from $C$
into $C$ and define $U(t,s)=\prod_{n=\nu(s)+1}^{\nu(t)}F_n$, the product representing composition of functions. Then
$U$ is an ES. If each $F_n$ is $M_n$-Lipschitz and the product $\prod_{n=1}^\infty M_n$ is bounded from above by $M$,
then $U$ is an \mles. For instance, if $F_n=(I+A_n)^{-1}$, where $\{A_n\}$ is a family of $m$-accretive operators on
$C$, then the piecewise constant interpolation of infinite products of resolvents defines a \ces.} \wbx\end{ejem}

\subsection{Almost-orbits and asymptotic equivalence}

\noi Let $V$ be an evolution system on $C$. A locally bounded trajectory of the form $t\mapsto V(t,s)x$ for $s$ and $x$ fixed is an {\em orbit} of $V$.  A locally bounded function $u:\R_+\to C$ is an {\em almost-orbit} of $V$ if
\begin{equation}
\limty{t}\|u(t+h)-V(t+h,t)u(t)\|=0\hskip20pt\hbox{uniformly in}\hskip10pt h\ge 0.
\end{equation}
Orbits and
almost-orbits have, essentially, the same asymptotic behavior.\\
Note the relation and difference with the notion of   asymptotic pseudotrajectories  where the convergence is uniform
on compact time intervals (\cite[Benaim and Hirsch]{BeH}, \cite[Benaim, Hofbauer and Sorin] {BHS}). The current concept
is more demanding but will allow for more precise results (convergence rather than properties on  the set of limit
points).

\begin{thm}\label{aseq}
Let  $V$ be an evolution system. For the weak topology assume either that $V$ is Lipschitz or $X$ is weakly complete
(weak Cauchy nets are weakly convergent\footnote{The spaces $\luno$ and $L^1$, as well as all reflexive Banach spaces,
have this property. It is not the case if $X$ contains $c_0$, though (see p. 88 in \cite[Li and Queff\'elec]{LiQ}).}).
If every orbit of $V$ converges weakly (resp. strongly), then so does every almost-orbit.
\end{thm}

\dem For the strong topology, let $u$ be an almost-orbit of $V$ and let $\eps>0$. By definition, there is $S>0$ such
that
$$\|u(t+h)-V(t+h,t)u(t)\|<\eps/4$$
for all $h\ge 0$ and $t\ge S$. Define $\zeta(S)=\limty{t}V(t,S)u(S)$ and choose $T>S$ such that
$\|V(t,S)u(S)-\zeta(S)\|<\eps/4$ for all $t\ge T$. Then
$$\|u(t+h)-\zeta(S)\| \le \|u(t+h)-V(t+h,S)u(S)\|+\|V(t+h,S)u(S)-\zeta(S)\|<\eps/2$$
for all $t\ge T$ and all $h\ge 0$. Thus  $\|u(t')-u(t)\| <\eps$ for all $t, t' \ge T$  so that  $u(t)$ is Cauchy and converges.\\
 It is clear that this argument is valid for the weak
topology if $X$ is weakly complete. If it is not the case but $V$ is $L$-Lipschitz, one defines
$\zeta(s)=w-\limty{t}V(t,s)u(s)$ and verifies that
$$\sup_{p\ge 0}\|\zeta(s+p)-\zeta(s)\|\le L\sup_{p\ge 0}\|u(s+p)-V(s+p,s)u(s)\|,$$
which tends to zero as $s\to\infty$ showing that $\zeta(s)$ converges strongly to some $\zeta$. Then one easily proves
that $u(t)$ converges weakly to $\zeta$ as $t\to\infty$.\bx

\noi A special case of Theorem \ref{aseq} was proved in \cite[Passty]{Pas}, when $V$ is defined by a semigroup of
contractions or if the almost-orbits are orbits of a semigroup of contractions.

\begin{thm}
Under the hypotheses of Theorem \ref{aseq}, the conclusion remains valid if the word {\em converges} is replaced
by {\em converges in average}.
\end{thm}

The proof of this result can be found in \cite[Alvarez and Peypouquet]{AlP2}.\\

\noi{\it Comments}\\
A similar result holds for {\em almost-convergence} (see \cite[Alvarez and Peypouquet]{AlP2}), a concept developed in
\cite[Lorentz]{Lor} that is stronger than convergence in average. It had been proved in \cite[Miyadera and
Kobayasi]{KoM} under supplementary assumptions: i) $V$ is defined by a strongly continuous semigroup of contractions;
ii) $\S\neq\emptyset$; and iii) for the weak topology, $X$ is weakly complete.\theend

\begin{lem}\label{lem:bounded_almost-orbit}
Let $U$ and $V$ be evolution systems and assume that for each $r>0$
$$\limty{t}\sup_{h\ge 0}\sup_{\|z\|\le r}\|U(t+h,t)z-V(t+h,t)z\|=0$$
then every bounded orbit of $V$ is an almost-orbit of $U$ and viceversa.
\end{lem}

\dem Let $v$ be an orbit of $V$ such that $\|v(t)\|\le r$ for all $t$. Then
\begin{eqnarray*}
\|v(t+h)-U(t+h,t)v(t)\|& = & \|V(t+h,t)v(t)-U(t+h,t)v(t)\|\\
& \le & \sup_{\|z\|\le r}\|U(t+h,t)z-V(t+h,t)z\|
\end{eqnarray*}
and so $v$ is an almost-orbit of $U$.\bx

\subsection{Continuous dynamics and discretizations}

\noi The following results explain why in most cases the systems defined in the preceding sections converge under the
same hypotheses. The proofs are considerably simplified if one assumes boundedness of the almost-orbits by virtue of
Lemma \ref{lem:bounded_almost-orbit}. We shall give them in this case  along with the references for more general
settings. The following proposition gathers results from \cite[Sugimoto and Koizumi]{KoS} and \cite[G\"uler]{Gul}.

\begin{prop}\label{prop:prox_vs_continuous} Let $A$ be a maximal monotone operator on $H$ and let  $U$ and $V$ be the evolution systems
defined by the differential inclusion  \eqref{continuous} and the proximal  algorithm \eqref{prox}, respectively. Assume one of the following
conditions holds:
\begin{itemize}
    \item [i)] $\{\lambda_n\}\in\ltwo\setminus\luno$; or
    \item [ii)] $A=\partial f$ and $\{\lambda_n\}\notin\luno$.
\end{itemize}
Then every orbit of $U$ is an almost-orbit of $V$ and viceversa.
\end{prop}

\dem Define $\nu(t)$ as in Example \ref{E:es2}. If $\{\lambda_n\}\in\ltwo\setminus\luno$, part $i)$ in Corollary
\ref{cor:prox_continuous} gives
$$\|U(t+s,t)z-V(t+s,t)z\|^2\le 3\|A^0z\|^2\sum_{n=\nu(t)}^\infty\lambda_n^2$$
and we conclude using Lemma \ref{lem:bounded_almost-orbit}. For unbounded almost-orbits, see \cite[Sugimoto and
Koizumi]{KoS}. If $A=\partial f$ and $\{\lambda_n\}\notin\luno$ the proof is highly technical and can be found in
\cite[G\"uler]{Gul}. It also relies on part $i)$ in Corollary \ref{cor:prox_continuous} but sharper estimations on
$\|A^0x_n\|$ and $\|A^0u(t)\|$ are needed. \bx

\begin{prop}
Let $T$ be nonexpansive, set $A=I-T$ and let $U$ and $W$ be the evolution systems defined by the differential inclusion  \eqref{continuous}
and Euler's discretization \eqref{euler}, respectively. Assume $\{\lambda_n\}\in\ltwo\setminus\luno$. Then every orbit of $U$ is an
almost-orbit of $W$ and viceversa.
\end{prop}

\dem The argument in the proof of part $i)$ in Proposition \ref{prop:prox_vs_continuous} can be applied here as well,
by virtue of inequality \eqref{Eq:Kobayashi_euler}.\bx

These properties allow for a better understanding of  similar asymptotic behavior of the continuous and discrete processes: in general for weak convergence in average (Section 4), for weak convergence in the case of demi-positive operators (Section 5) and for strong convergence under addtitional geometrical hypotheses (Section 6).

\subsection{Quasi-autonomous systems}

One of the advantage of this approach through almost-orbits  is that it extends to non-autonomous systems.

\subsubsection{Continuous dynamics}

Recall that the solutions of the differential inclusion \eqref{continuous} define an evolution system $U$ as in Example
\ref{E:es1}. Let us consider quasi-autonomous versions of \eqref{continuous}, namely
\begin{equation}\label{Eq:continuous-quasi}
-\dot v(t)\in Av(t)+\f(t)
\end{equation}
and
\begin{equation}\label{Eq:continuous-Tik}
-\dot v(t)\in Av(t)+\eps(t)v(t).
\end{equation}

\begin{prop}\label{perl1}
If $\f\in L^1(0,\infty;X)$, then every function satisfying \eqref{Eq:continuous-quasi} is an almost-orbit of $U$. The
same holds for every function satisfying \eqref{Eq:continuous-Tik} provided $\eps\in L^1(0,\infty;\R)$.
\end{prop}

\dem For the first part we follow \cite[Miyadera and Kobayasi]{KoM}. If $v$ satisfies \eqref{Eq:continuous-quasi} and
$t\ge 0$ we have
$$
\|v(t+s)-U(t+s,t)v(t)\|^2\le 2 \int_0^s\|\f(t+\tau)\|\|v(t+\tau)-U(t+\tau,t)v(t)\|d\tau
$$
and so
$$
\|v(t+s)-U(t+s,t)v(t)\|\le\int_0^s\|\f(t+\tau)\|d\tau\le\int_t^\infty\|\f(\tau)\|d\tau.
$$
On the other hand, let $v$ satisfy \eqref{Eq:continuous-Tik}. Fix $t$ and consider as above
$\psi(s)=\frac{1}{2}\|U(t+s,t)v(t)-v(t+s)\|^2$. Using
$\langle\zeta,\zeta-\xi\rangle\ge-\frac{1}{4}\|\xi\|^2$ for all $\zeta,\xi\in H$, we deduce
$\dot\psi(s)\le\frac{1}{4}|\eps(t+s)|\|U(t+s,t)v(t)\|^2$ for almost every $s>0$. Integrating from $0$ to $s$ and observing
that $\psi(0)=0$ we obtain
$$\|U(t+s,t)v(t)-v(t+s)\|^2\le\frac{1}{4}\int_t^{t+s}|\eps(\tau)|\ \|U(t+\tau,t)v(t)\|^2\ d\tau\le\frac{M}{4}\int_t^{\infty}|\eps(\tau)|\ d\tau$$
if $v$ is bounded. \bx

{\it Comments}\\
In \cite[Alvarez]{Alv}, the author studies the problem
\begin{equation}\label{osc}
u''(t)+\gamma u'(t)+\nabla\Phi(u(t))=0,
\end{equation}
where $\Phi$ is a $\C^1$ convex function. He proves that if $\argmin(\Phi)\neq\emptyset$, then each solution $u(t)$
converges weakly to a minimizer of $\Phi$ as $t\to\infty$ and gives conditions for strong convergence. Later, in
\cite[Attouch and Czarnecki]{ACz} the authors establish, among other results,  that if $\eps\in L^1$ the solutions of
\begin{equation}\label{osceps}
u''(t)+\gamma u'(t)+\nabla\Phi(u(t))+\eps(t)u(t)=0.
\end{equation}
also converge weakly to minimizers of $\Phi$. It turns out (see \cite[Alvarez and Peypouquet]{AlP3}) that under this
condition ($\eps\in L^1$)  the solutions of \eqref{osceps} are
almost-orbits of the evolution system defined by \eqref{osc}. \\
This is an alternative way to prove the cited result from \cite[Attouch and Czarnecki]{ACz} and  it shows that these
tools building on almost-orbits to classify the asymptotic behavior through  equivalence classes (continuous
trajectories, proximal or Euler approximations, Tykhonov regularization, perturbations) can be applied to second-order
systems as well.\theend

\subsubsection{Proximal sequences}

In a similar fashion one can prove any interpolation of a sequence $\{y_n\}$ satisfying
\begin{equation}\label{Eq:prox-quasi}
y_{n-1}-y_n\in \lambda_nAy_n+\phi_n
\end{equation}
or
\begin{equation}\label{Eq:prox-Tik}
y_{n-1}-y_n\in \lambda_nAy_n+\epsilon_ny_n
\end{equation}
is an almost-orbit of the evolution system $U$ defined by the proximal scheme \eqref{prox} as in Example \ref{E:es2}
provided $\{\phi_n\}\in\luno(\N;X)$ and $\{\epsilon_n\}\in\luno(\N;\R_+)$, respectively.\\

For additional applications and examples see \cite[Alvarez and Peypouquet]{AlP3}.\\


\section{Concluding remarks}

It is useful  to observe that there are two aspects related to  the ideas of asymptotic equivalence discussed in the last
section. In the first place, one can obtain sufficient conditions for a perturbed, regularized or discretized system to
have the same asymptotic properties as the original one.  The issue here is in terms of stability or regularity or computational purposes. On the other hand, if a given dynamics  does not have some
desirable asymptotic behavior,  one can introduce pertubation in order to generate  orbits having better properties. In this case,
the tools  of asymptotic equivalence give necessary condition for a perturbation to be effective.\\

Observe that the trajectories defined by \eqref{continuous} only converge weakly in average. Even in the case where
$A=\partial f$, convergence is still weak and the limit depends on the initial point. One can get a better asymptotic
behavior by forcing the system to stabilize in the direction of the origin. More precisely, consider a piecewise
absolutely continuous function $\eps:\R_+\to\R_+$ such that $\limty{t}\eps(t)=0$. If $\eps\in L^1(0,\infty;\R_+)$ the
system defined by \eqref{Eq:continuous-Tik} will have the same asymptotic behavior as \eqref{continuous} by Proposition
 \ref{perl1}. If we expect
the regularized system to have better properties we must consider $\eps\notin L^1(0,\infty;\R_+)$. The following result
is from \cite[Cominetti, Peypouquet and Sorin]{CPS}:

\begin{prop}
Suppose $v:\R_+\to H$ satisfies
$$-\dot v(t)\in Av(t)+\eps(t)v(t).$$
with $\eps\notin L^1(0,\infty;\R_+)$. Assume further that $A=\partial f$ or $\int_0^\infty|\dot\eps(t)|\ dt<\infty$ (finite total variation). Then
$\limty{t}v(t)=P_{\S}0$.
\end{prop}

\noi Special cases of the preceding result had been proved earlier in \cite[Browder]{Bro}, \cite[Reich]{Rei2} and
\cite[Attouch and Cominetti]{AtC}.
A similar result for the second order appears in  \cite[Attouch and Czarnecki]{ACz}. \\

Also, as we mentioned before, the trajectories defined by \eqref{continuous} need not be weakly convergent. If one applies
the proximal point algorithm with stepsizes $\lambda_n\in\ltwo$, by Proposition \ref{prop:prox_vs_continuous}, the
corresponding system will have the same asymptotic properties. In other words, the approximation is too good: ``the discrete approximation mirrors the behavior of the differential equation {\it too} well" \cite[Bruck, p. 29]{BruS}. If one
wishes to get a better (or different) behavior, it is necessary to consider $\lambda_n\notin\ltwo$. This turns out to
be fruitful because, in that case Theorem \ref{weakprox} guarantees weak convergence even when the operator is not
demipositive (see also Example 3 in Section 6).\\

{\bf Acknowledgments}\\
The authors wants to thank Roberto Cominetti for his help during the preparation of this work. In particular his
(unpublished) notes on ``Evolution equations and monotone maps" were very helpful.\\
S. Sorin acknowledges support from grant ANR-08-BLAN-0294-01 (France).


\end{document}